\documentclass[12pt]{article}
\usepackage{amsfonts,amsmath,appendix,color,amssymb,graphicx}
\DeclareGraphicsExtensions{.pdf}

\makeatletter
\@addtoreset{equation}{section}
\makeatother
\newtheorem{theorem}{Theorem}[section]

\newtheorem{remark}{Remark}[section]
\newtheorem{proposition}{Proposition}[section]
\newtheorem{hypothesis}{Hypothesis}[section]
\newtheorem{lemma}{Lemma}[section]
\newtheorem{corollary}{Corollary}[section]
\newtheorem{definition}{Definition}[section]
\def\br{\begin{remark}\rm\small}
\def\er{\end{remark}}
\def\bt{\begin{theorem}}
\def\et{\end{theorem}}
\def\bd{\begin{definition}}
\def\ed{\end{definition}}
\def\bp{\begin{proposition}}
\def\ep{\end{proposition}}
\def\bl{\begin{lemma}}
\def\el{\end{lemma}}
\def\bc{\begin{corollary}}
\def\ec{\end{corollary}}
\def\beaq{\begin{eqnarray}}
\def\eeaq{\end{eqnarray}}

\newcommand{\nn}{\nonumber}
\newcommand{\dd}{\mathrm{d}}
\newcommand{\beq}{\begin{equation}}
\newcommand{\eeq}{\end{equation}}
\newcommand{\bea}{\begin{eqnarray}}
\newcommand{\eea}{\end{eqnarray}}
\newcommand{\n}{\mathfrak{n}}
\newcommand{\ra}{\rightarrow}
\newcommand{\e}{\epsilon}

\newcommand{\p}{{\parallel}}

\newcommand{\Tr}{{\,\rm Tr}\:}

\newcommand{\Res}{\mathop{\,\rm Res\,}}

\definecolor{rouge}{rgb}{0.84,0.18,0.07}
\definecolor{bleu}{rgb}{0.22,0.41,0.74}
\definecolor{vertf}{rgb}{0.08,0.46,0.07}
\usepackage[pdftex]{hyperref}
\hypersetup{colorlinks,urlcolor=vertf,citecolor=rouge,linkcolor=bleu,filecolor=black}

\begin{document}

\sloppy

\thispagestyle{empty}
\addtolength{\baselineskip}{0.20\baselineskip}
\begin{center}
\vspace{26pt}
{\large \bf {Asymptotic expansion of $\beta$ matrix models
in the one-cut regime}}
\end{center}

\vspace{26pt}

\begin{center}
{\sl Ga\"etan Borot}\hspace*{0.05cm}\footnote{\href{mailto:gaetan.borot@cea.fr}{gaetan.borot@cea.fr}},
{\sl Alice Guionnet}\hspace*{0.05cm}\footnote{\href{mailto:alice.guionnet@ens-lyon.fr}{alice.guionnet@ens-lyon.fr}}

\vspace{6pt}
$\,^{1}$ Institut de Physique Th\'{e}orique de Saclay,\\
91191 Gif-sur-Yvette Cedex, France, \\
and Section de Math\'ematiques, Universit\'e de Gen\`eve \\
2-4 rue du Li\`evre, 1211 Gen\`eve 4, Suisse. \\

\vspace{0.2cm}

$\,^{2}$ UMPA, CNRS  UMR 5669, ENS Lyon, \\
46 all\'ee d'Italie,
69007 Lyon, France.\\

\end{center}

\vspace{20pt}
\begin{center}
{\bf Abstract}
\end{center}

%

\vspace{0.5cm}

We prove the existence of a $1/N$ expansion to all orders in $\beta$ matrix models with a confining, offcritical potential corresponding to an equilibrium measure with a connected support. Thus, the coefficients of the expansion can be obtained recursively by the "topological recursion" derived in \cite{CE06}. Our method relies on the combination of a priori bounds on the correlators and the study of Schwinger-Dyson equations, thanks to the uses of classical complex analysis techniques. These a priori bounds can be derived following \cite{Boutet,Johansson, Shch2}
or for strictly convex potentials by using  concentration of measure \cite[Section 2.3]{AGZ}. Doing so, we extend the strategy of \cite{GMS}, from the hermitian models ($\beta = 2$) and perturbative potentials, to general $\beta$ models. The existence of the first correction in $1/N$ was considered in \cite{Johansson} and more recently in \cite{Shch2}. Here, by taking similar hypotheses, we extend the result to all orders in $1/N$.

\vspace{0.5cm}

\newpage


\vspace{26pt}
\pagestyle{plain}
\setcounter{page}{1}
\addtocounter{footnote}{-2}


\section{Introduction}
\label{sec:intro}

\subsection{Definitions}

We consider the probability measure $\mu_{N,\beta}^V$ on $\mathbb R^N$
given by:
\beq
\label{eqmes}\mathrm{d}\mu_{N,\beta}^{V;[b_-,b_+]}(\lambda) = \frac{1}{Z_{N,\beta}^{V;[b_-,b_+]}}\prod_{i = 1}^N \mathrm{d}\lambda_i\,e^{-\frac{N\beta}{2} \,V(\lambda_i)}\,\mathbf{1}_{[b_-,b_+]}(\lambda_i)\,\prod_{1 \leq i < j \leq N} |\lambda_i - \lambda_j|^{\beta}.
\eeq
 $[b_-,b_+]$ is an interval of the real line, $-\infty \leq b_-<b_+ \leq +\infty$, and $\beta$ is a positive number. For $\beta = 1,2,4$, this is the measure
induced on the eigenvalues of $\Phi$ by the probability measure
$\mathrm{d}\Phi\,e^{-\frac{N\beta}{2}\Tr V(\Phi)}$ on a vector space
$\mathcal{E}_{N,\beta}$ of $N \times N$ matrices.
$\mathcal{E}_{N,1}$ (resp. $\mathcal{E}_{N,2}$, and
$\mathcal{E}_{N,4}$) is the space of real symmetric (resp. hermitian,
and quaternionic self-dual) matrices \cite{Mehtabook}. For general
$\beta > 0$, when $V$ is quadratic (Hermite weight), or log + linear (Laguerre weight), Dumitriu and Edelman have found \cite{DE02} a measure on a set of tridiagonal matrices which induces the measure $\dd\mu_{N,\beta}^{V;\mathbb R}$ on eigenvalues. This can be generalized to any even polynomial potential \cite{Riderco}. This was subsequently exploited to study these particular $\beta$ matrix models by Ram\`{i}rez, Rider and Vir\'{a}g \cite{RRV}.  Though, there is no known plain random matrix whose spectrum is distributed according to $\dd\mu_{N,\beta}^{V;\mathbb{R}}$
for general $\beta$ and $V$, we still speak of a "matrix model", and we call $\lambda_i$ the eigenvalues.

We define the unnormalized empirical measure $M_N$ of the eigenvalues
given by $$M_N=\sum_{i=1}^N \delta_{\lambda_i}$$ and their Cauchy-Stieltjes transform, which are generating series for the moments of $M_N$. In fact, we prefer to work with the generating series of the cumulants. They are also called "correlators", and are
defined for $x_1,\ldots,x_n\in\mathbb C\backslash \mathbb [b_-,b_+]$ by
\bea
W_n^{V;[b_-,b_+]}(x_1,\ldots, x_n)& = & \mu_{N,\beta}^{V;[b_-,b_+]}\Big[\Big( \int \frac{\dd M_N(\xi_1)}{x_1-\xi_1}\cdots\int\frac{\dd M_N(\xi_n)}{x_n - \xi_n}\Big)_c\Big] \nn \\
& = & \partial_{\epsilon_1}\cdots\partial_{\epsilon_n}\Big(\ln Z_{N,\beta}^{V-\frac{2}{\beta N}\sum_i \frac{\epsilon_i}{x_i- \bullet};[b_-,b_+]}\Big)\Big|_{\epsilon_i=0}. \nn
\eea
In particular, we have
{\small\bea
W_1^{V;[b_-,b_+]}(x) & = & \mu^{V;[b_-,b_+]}_{N,\beta}\Big[\int \frac{\dd M_N(\xi)}{x-\xi} \Big],\nn \\
W_2^{V;[b_-,b_+]}(x_1,x_2) & = & \mu^{V;[b_-,b_+]}_{N,\beta}\Big[\iint \frac{\dd M_N(\xi)}{x_1 - \xi}\,\frac{\dd M_N(\eta)}{x_2 - \eta}\Big]\nn\\
&&\quad - \mu^{V;[b_-,b_+]}_{N,\beta}\Big[\int \frac{\dd M_N(\xi)}{x_1 - \xi}\Big]\,\mu_{N,\beta}^{V;[b_-,b_+]}\Big[\int \frac{\dd M_N(\eta)}{x_2 - \eta}\Big] .\nn
\eea}
$\!\!\!$When there is no confusion, we may omit to write the dependence in $V$ and $[b_-,b_+]$ in the exponent.

It is well-known, see \cite{Johansson} or the textbooks \cite[Theorem6]{MR1677884} or \cite[Theorem 2.6.1 and Corollary
2.6.3]{AGZ}, that
\begin{theorem}
\label{th:1} Assume that $V\,:\, [b_-,b_+] \rightarrow \mathbb{R}$ is a continuous function, and if $b_{\tau} = \tau\infty$ is infinite, assume that:
\beq
\liminf_{x \rightarrow \tau\infty} \frac{V(x)}{2\ln|x|} > 1 .\nn
\eeq
If $V$ depends on $N$, assume also that $V \rightarrow V^{\{0\}}$ in the space of continuous function over $[b_-,b_+]$ for the sup norm.
Then, the normalized empirical measure $L_N=N^{-1}\,M_N$ converges almost surely and in expectation
towards the unique probability measure $\mu_{\mathrm{eq}} := \mu_{\mathrm{eq}}^{V;[b_-,b_+]}$ on $[b_-,b_+]$ which minimizes:
\beq
\mathcal{E}[\mu] = \int \dd\mu(\xi)V^{\{0\}}(\xi)-\iint \dd\mu(\xi)\dd\mu(\eta)\ln|\xi-\eta|. \nn
\eeq
Moreover, $\mu_{\mathrm{eq}}$ has compact support. It is characterized by the existence of a constant $C$ such that:
\beq
\label{ina}\left\{\begin{array}{ll} \forall x \in [b_-,b_+] & \,\, 2\int_{b_-}^{b_+}\dd \mu_{\mathrm{eq}}(\xi)\ln|x - \xi| - V^{\{0\}}(x) \leq C, \\
\mu_{\mathrm{eq}}\,\,\mathrm{almost}\,\,\mathrm{surely} & \,\, 2\int_{b_-}^{b_+} \dd \mu_{\mathrm{eq}}(\xi)\ln|x - \xi| - V^{\{0\}}(x) = C. \end{array}\right.
\eeq
\end{theorem}
In particular, for any $x\in\mathbb{C}\setminus[b_-,b_+]$, we have
\begin{equation}
\lim_{N\rightarrow\infty} \frac{1}{N}\,W_1(x) = \int\frac{\dd \mu_{\mathrm{eq}}(\xi)}{x-\xi} := W_{1}^{\{-1\}}(x), \nn
\end{equation}
and the convergence is uniform in any compact of $\mathbb{C}\setminus[b_-,b_+]$.

\subsection{Main results}

Our goal is to prove an asymptotic expansion in powers of $1/N$ when $N \rightarrow \infty$  for the partition function $Z_{N,\beta}^{V;[b_-,b_+]}$ and the correlators $W_n^{V;[b_-,b_+]}(x_1,\ldots,x_n)$. This is not always expected. In particular it is false when the support of $\mu_{\mathrm{eq}}^{V;[b_-,b_+]}$, the limiting eigenvalue distribution, is not connected: corrections to the leading order feature a quasi periodic behavior with $N$ (see \cite{Ecv} for a general heuristic argument). Our proof uses a priori bounds on the correlators, and what we really need is to establish that $W_n \in O(1)$ for $n \geq 2$. We shall prove this condition either
based on a result of Boutet de Monvel, Pastur and Shcherbina \cite{Boutet} (also used recently in the context of $\beta$ ensembles by Kriecherbauer and Shcherbina \cite{Shch2}), or under the additional assumption that $V$ is strictly convex. In the convex setting, such a priori bounds can be derived from concentration of measures properties, in which case our article is self-contained. Our basic assumptions and main results are:

\begin{hypothesis}
\label{hypoge}
\begin{itemize}
\item[] \phantom{sss}
\item[$\bullet$] (Regularity) $V\,:\,[b_-,b_+] \rightarrow \mathbb{R}$ is continuous, and if $V$ depends on $N$, it has a limit $V^{\{0\}}$ in the space of continuous functions over $[b_-,b_+]$ for the sup norm.
\item[$\bullet$] (Confinement) If $b_{\tau} = \tau\infty$, $\liminf_{x \rightarrow \tau\infty} \frac{V(x)}{2\ln|x|} > 1$.
\item[$\bullet$] (One-cut regime) The support of $\mu_{\mathrm{eq}}^{V;[b_-,b_+]}$ consists in a unique interval $[\alpha_{-},\alpha_+] \subseteq [b_-,b_+]$.
\item[$\bullet$] (Control of large deviations) The function $x\in [b_-,b_+]\backslash ]\alpha_-,\alpha_+[ \mapsto \frac{1}{2}V(x)-\int\ln|x-\xi|d\mu_{\rm eq}(\xi)$ achieves its minimum value at $\alpha_-$ and $\alpha_+$ only.
\item[$\bullet$] (Offcriticality)
$S(x) > 0$ whenever $x \in [\alpha_-,\alpha_+]$, where:
\beq
S(x) = \pi\,\frac{\dd \mu_{\mathrm{eq}}}{\dd x}\,\sqrt{\left|\frac{\prod_{\tau' \in \mathsf{Hard}} (x - \alpha_{\tau'})}{\prod_{\tau \in \mathsf{Soft}} (x - \alpha_{\tau})}\right|}, \nn
\eeq
where $\tau\in \mathsf{Hard}$ (resp. $\tau\in \mathsf{Soft}$)
iff $b_\tau=\alpha_\tau$ (resp. $\tau (b_\tau-\alpha_\tau)>0$).
\item[$\bullet$] (Analyticity) $V$ can be extended as a holomorphic function in some open neighborhood of $[\alpha_-,\alpha_+]$.
\item[$\bullet$] $V$ has a $1/N$ expansion in this neighborhood, in the sense of Hyp.~\ref{hypVN}.
\end{itemize}
\end{hypothesis}
Notice that the "one-cut regime", "offcriticality" and "control of large deviations" assumptions automatically hold when $V$ is strictly convex
(see \cite[Proposition 3.1]{Johansson}, which extends easily to analytic
functions instead of polynomials).

\begin{proposition}
\label{prop1}
Assume Hyp.~\ref{hypoge}. Then, $W_n^{V;[b_-,b_+]}$ admits an asymptotic expansion when $N \rightarrow \infty$:
\beq
W_n^{V;[b_-,b_+]}(x_1,\ldots,x_n) = \sum_{k \geq n - 2} N^{-k}\,W_n^{V;\{k\}}(x_1,\ldots,x_n), \nn
\eeq
which has the precise meaning that, for all $K \geq n - 2$:
\beq
W_n^{V;[b_-,b_+]}(x_1,\ldots,x_n) = \sum_{k = n - 2}^{K} N^{-k}\,W_n^{V;\{k\}}(x_1,\ldots,x_n) + o(N^{-K}). \nn
\eeq
The $o(N^{-K})$ is uniform for $x_1,\ldots,x_n$ in any compact of $(\mathbb{C}\setminus[b_-,b_+])^n$, but not uniform in $n$ and $K$. Moreover, if $(b_{\tau} - \alpha_{\tau})\tau > 0$ (meaning that $\alpha_{\tau}$ is a soft edge) and Hypotheses~\ref{hypoge} hold, the functions $W_n^{V;\{k\}}$ are independent of $b_{\tau}$.
\end{proposition}

\begin{proposition}
\label{prop0}
Assume Hyp.~\ref{hypoge}, and $b_- < \alpha_- < \alpha_+ < b_+$ (all edges are soft). Then, $Z_{N,\beta}^{V;[b_-,b_+]}$ admits an asymptotic expansion when $N \rightarrow \infty$:
\beq
\label{eq:expxp} Z_{N,\beta}^{V;[b_-,b_+]} = Z_{N,\mathrm{G}\beta\mathrm{E}} \left(\frac{\alpha_+ - \alpha_-}{4}\right)^{N + \beta\frac{N(N - 1)}{2}}\,\exp\left(\sum_{k \geq -2} N^{-k}\,F^{V;\{k\}}_{\beta}\right).
\eeq
In other words, for any $K \geq -2$:
\beq
Z_{N,\beta}^{V;[b_-,b_+]} = Z_{N,\mathrm{G}\beta\mathrm{E}} \left(\frac{\alpha_+ - \alpha_-}{4}\right)^{N + \beta\frac{N(N - 1)}{2}}\,\exp\left(\sum_{k = -2}^{K} N^{-k}\,F^{V;\{k\}}_{\beta} + o(N^{-K})\right). \nn
\eeq
Moreover, the coefficients $F^{V;\{k\}}_{\beta}$ are independent of $b_-$ and $b_+$ chosen such that $b_- < \alpha_- < \alpha_+ < b_+$ and Hypotheses~\ref{hypoge} hold.
\end{proposition}

$Z_{N,\mathrm{G}\beta\mathrm{E}}$ is the partition function of the Gaussian $\beta$ ensemble, defined by the quadratic potential $V_{\mathrm{G}}(x) = \frac{x^2}{2}$. It is given by a Selberg integral \cite{Sel} (also in \cite{Mehtabook}):
\beq
\label{eq:Selb} Z_{N,\mathrm{G}\beta\mathrm{E}} = (2\pi)^{N/2}\,\big(N\beta/2\big)^{-\beta N^2/4 + (\beta/4 - 1/2)N}\,\frac{\prod_{j = 1}^N \Gamma\big(1 + j\beta/2\big)}{\Gamma\big(1 + \beta/2)^N}.
\eeq
For hard edges (i.e. $b_- = \alpha_-$ or $b_+ = \alpha_+$), one may still interpolate between $Z_{N,\beta}^{V;[b_-,b_+]}$ and a Gaussian $\beta$ ensemble restricted to some interval (Corollary~\ref{interpol}), but the partition function of the latter is not a Selberg integral and thus not known in closed form.

\subsubsection{Commentary}

When $V$ does not depend on $N$ and $\beta$, $W_n^{V;\{k\}}$ has a very simple dependence in $\beta$:
\beq
\label{nwn}W_n^{V;\{k\}} = \sum_{g = 0}^{\lfloor (k - n + 2)/2 \rfloor} \left(\frac{\beta}{2}\right)^{1 - g - n}\,\Big(1 - \frac{2}{\beta}\Big)^{k + 2 - 2g - n}\,\mathcal{W}_n^{V;(g;k + 2 - 2g - n)},
\eeq
and likewise:
\beq
\label{nwnq}F_{\beta}^{V;{\{k\}}} =  \sum_{g = 0}^{\lfloor k/2 \rfloor + 1} \left(\frac{\beta}{2}\right)^{1 - g}\,\Big(1 - \frac{2}{\beta}\Big)^{k + 2 - 2g}\,\mathcal{F}^{V;(g;k + 2 - 2g)}.
\eeq
Assuming existence of the $1/N$ expansion, or at the level of formal matrix integrals, the recursive computation of the ${\mathcal W}_n^{V;(g;l)}$ and $\mathcal{F}^{V;(g;l)}$ was developed by Chekhov and Eynard in \cite{CE06}. For $\beta = 2$, it is well-known that Eqn.~\ref{eq:expxp} is an expansion in even powers of $N$, i.e. $F^{V;\{2k + 1\}}_{\beta = 2} = 0$. Such a result goes back to the so-called topological expansion of t'Hooft, shown in the context of matrix models by Br\'{e}zin, Itzykson, Parisi and Zuber \cite{BIPZ}. Indeed when $\beta = 2$, the sum in Eqn.~\ref{nwn} has only one term, namely $k = 2g - 2 + n$, which is present only when $k = n\,\,\mathrm{mod}\,\,2$, and likewise for Eqn.~\ref{nwnq} which can be considered as the case $n = 0$.

At the asymptotic level, the case $\beta=2$ was tackled in \cite{APS01}. For $\beta = 1,2,4$, the partition function and the correlators can be computed with the help of orthogonal polynomials \cite{Mehtabook}. These are solutions of a Riemann-Hilbert problem \cite{FIK92}, for which the large $N$ asymptotics have been intensively studied \cite{BI99,DKMcLVZ1,DKMcLVZ2,DKMcLVZ3} with the steepest descent method introduced in \cite{DZ}. As a consequence, Ercolani and
 McLaughlin \cite{ErcMcL} were able to prove the existence of a $1/N^2$ expansion of $\ln Z_{N,\beta = 2}^V$. However, the topological expansions in the cases $\beta = 1$ and $4$ are technically more involved in this framework, and have resisted to analysis up to now. 

Integrability properties of $\beta$ matrix models are unraveled for general $\beta > 0$, in particular there is no known orthogonal polynomials techniques to evaluate the partition function $Z_{N,\beta}^V$ and the correlation functions $W_n(x_1,\ldots,x_n)$. Yet, it is always possible to study the Cauchy-Stieltjes transform of the empirical measure of the eigenvalues and the  "loop equations", also called Schwinger-Dyson equations or Pastur equations \cite{Pastur72}, that govern its expectations and cumulants. Thanks to the rough bounds for $W_1^{V}$ and $W_2^{V}$ established in \cite{Boutet}, Johansson \cite{Johansson} proved a central limit theorem and obtained the first correction to $W_1$ when $V$ is an even polynomial satisfying Hyp.~\ref{hypoge}. This was also the subject of a recent work by Kriecherbauer and Shcherbina \cite{Shch2}, with Hyp.~\ref{hypoge} only. These authors have obtained in particular the expansion of $\ln Z_{N,\beta}^{V}$ up to a $O(1)$ when $N \rightarrow \infty$ (see their Theorem 2).

The determination of $W_1^{V;\{-1\}}$ \cite{Wig,BIPZ,BAG} and $W_2^{V;\{0\}}$ \cite{AM90,Been2,Johansson} has been known for long, in $\beta$ ensembles or many other matrix models.
It was also observed long ago \cite{AM90} that, if a $1/N$ expansion is assumed to exist, the loop equations turn into a system of recursive linear equations determining fully the decaying orders. To solve it, one just has to invert a linear operator $\mathcal{K}$. Recursiveness is a consequence of the assumption or the fact that $W_n \in O(N^{2 - n})$, which allows the determination of the leading order of $W_{n - 1}$ without knowledge of $W_n$ (for $n \geq 3)$. These techniques found their origin in \cite{AM90,ACM,ACKM,ACKMe} and culminated with the formalism of the "topological recursion" of \cite{E1MM,EOFg} for $\beta = 2$, and \cite{CE06} for any fixed $\beta > 0$.

In this article, we observe that $\mathcal{K}^{-1}$ is a continuous operator on some appropriate space of analytic functions. Combining with the a priori control on correlators which dates back to \cite{Boutet}, we prove the existence of the full expansion.

For strictly convex potentials, concentration inequalities also provide rough bounds on the correlators and therefore allow us to give self-contained proofs,independent from \cite{Boutet}-\cite{Shch2}. In this framework, loop equations were used in \cite{GMS} to establish the asymptotic expansion of models of several hermitian random matrices ($\beta = 2$) with strictly convex interactions. Maurel-Segala \cite{MSp} also studied models of several symmetric random matrices ($\beta = 1$) with strictly convex interactions. In order to prove the asymptotic expansion, the main step of \cite{GMS} was to show that some operator on non-commutative polynomials could be inverted, with bounded appropriate norm, and this was only done in a perturbative regime. Here, thanks to complex analysis, the potential need not be a small perturbation of the quadratic potential.

Our techniques could also be applied to other matrix models. For instance, the convergent $\beta$, $\mathcal{O}(\mathfrak{n})$ matrix model:
\beq
\dd\mu_{N,\beta,\mathcal{O}(\mathfrak{n})}^{V;\mathbb{R}_+}(\lambda) = \frac{1}{Z_{N,\beta,\mathcal{O}(\mathfrak{n})}^{V;\mathbb{R}_+}}\,\prod_{i = 1}^{N} \dd \lambda_i\,e^{-\frac{N\beta}{2}\,V(\lambda_i)}\,\frac{\prod_{1 \leq i < j \leq N}
 |\lambda_i - \lambda_j|^{\beta}}{\prod_{1 \leq i,j \leq N} (\lambda_i + \lambda_j)^{\mathfrak{n}/2}}. \nn
\eeq
An important point is that the corresponding quadratic functional:
\beq
\mathcal{E}[\rho] = \iint_{\mathbb{R}_+^2} \dd\rho(\xi)\dd\rho(\eta)\Big[-\frac{\beta}{2}\,\ln|\xi - \eta| + \frac{\mathfrak{n}}{2}\,\ln|\xi + \eta|\Big] + \frac{\beta}{2}\int_{\mathbb{R}_+}\dd\rho(\xi)\,V^{\{0\}}(\xi) \nn
\eeq
is strictly convex in the regime $|\n| < \beta$, therefore ensuring uniqueness of its minimizer. Besides, the analytic tools for the recursive determination of the one-cut solution to the loop equations of the $\mathcal{O}(\mathfrak{n})$ model in this regime were clarified in \cite{BEOn}. The existence of a full $1/N$ expansion for convergent $\mathcal{O}(\mathfrak{n})$ matrix models under Hyp.~\ref{hypoge} could probably be established by following the lines we are presenting for the $\beta$ matrix models.

An open challenge, which would be interesting for further applications, is to extend the same reasoning to complex measures, i.e. to Eqn.~\ref{eqmes} where the potential $V$ is complex-valued, or/and where the eigenvalues are integrated on a given path in the complex plane.

\subsection*{Outline of the article}

We first study in Section~\ref{weak} the weak dependence in the bounds of integration under weak assumptions on $V$. In particular, we may trade the initial interval $[b_-,b_+]$ for a finite interval $[a_-,a_+]$. We then write in Section~\ref{secloop} the corresponding loop equations for the correlators. Section~\ref{sec:rec1} is devoted to the proof of the asymptotic expansion of the correlators with slightly stronger assumptions (Prop.~\ref{recerror}). They are weakened in Section~\ref{final} to complete the proof of our main results for the correlators (Prop.~\ref{prop1}) and the free energy (Prop.~\ref{prop0}). We also remind how early steps of our proof imply the central limit theorem of Johansson (Prop.~\ref{propclt}).

\section{Weak dependence on the soft edges}
\label{weak}
In this section we show that the partition function and the correlators
depend very weakly on the boundary points
 of the interval of integration $[b_-,b_+]$
if they are soft, i.e. do not coincide with the boundary points of the
 support  $[\alpha_-,\alpha_+]$
of the equilibrium measure. We show more precisely that
this dependence yields only exponentially small corrections, by deriving a large deviation principle
for the law of the extreme eigenvalues.  This point was already  studied in \cite[section 2.6.2]{AGZ}
under a technical assumption \cite[Assumption 2.6.5]{AGZ} that we replace here by assuming that the rate function of our large deviation principle vanishes only at
$\alpha_-$ and $\alpha_+$.  This result is not new in essence and not specific to the one-cut regime, see for instance \cite[Proposition 2]{APS01} where it is proved with the extra assumption that $V$ has bounded second derivatives in a neighborhood of $\mathrm{supp}\,\mu^{\mathrm{eq}}$, or \cite[Proposition 11.1.4]{PSbook} where it is proved with the extra assumption that $V$ satisfies a Lipschitz condition in $[b_-,b_+]$. 

\subsection{Large deviation principle for the extreme eigenvalues}

Under the assumptions of Theorem~\ref{th:1} on an interval $[b_-,b_+]$, we can define:
\beq
\mathcal{J}^{V;[b_-,b_+]}(x)  = \frac{V(x)}{2} - \int_{b_-}^{b_+} \dd \mu_{\mathrm{eq}}^{V;[b_-,b_+]}(\xi)\,\ln|x - \xi| \nn
\eeq
when $x \in [b_-,b_+]$, and $+\infty$ otherwise. Suppose that $[b_-,b_+] \neq [\alpha_-,\alpha_+]$, and set:
\beq
\widetilde{\mathcal{J}}^{V;[b_-,b_+]}(x) = \mathcal{J}^{V;[b_-,b_+]}(x) - \inf_{\xi \in [b_-,b_+]} \mathcal{J}^{V;[b_-,b_+]}(\xi). \nn
\eeq
We define also $\widetilde{\mathcal{J}}_{\mathrm{max}}^{V;[b_-,b_+]}(x)$ (resp. $\widetilde{\mathcal{J}}^{V;[b_-,b_+]}_{\mathrm{min}}(x)$) which is equal to $\widetilde{\mathcal{J}}^{V;[b_-,b_+]}(x)$, except when $x \in ]-\infty,\alpha_-]$ (resp. $[\alpha_+,+\infty[$) where we set its value to  $+\infty$.

\begin{proposition}
\label{uuu3} Let $V\::\:[b_-,b_+] \rightarrow \mathbb{R}$ be a continuous function, and if $b_{\tau} = \tau \infty$, assume that:
\beq\label{assss}
\liminf_{x \rightarrow \tau\infty} \frac{V(x)}{2\ln|x|} > 1.
\eeq
Assume that $\widetilde{\mathcal{J}}^{V;[b_-,b_+]}$ does not vanish outside $[\alpha_-,\alpha_+]$. Then:
 \begin{itemize}
\item[$(i)$] $\beta\widetilde{\mathcal{J}}^{V;[b_-,b_+]}_{\rm max}$ (resp.
$\beta\widetilde{\mathcal{J}}^{V;[b_-,b_+]}_{\rm min}$) is
a  good rate function on $[b_-,b_+]$, which vanishes at $\alpha_{+}$ (resp. $\alpha_-$).
\item[$(ii)$] The law of $\lambda_{\mathrm{max}}$ (resp. $\lambda_{\mathrm{min}}$) under $\mu^{V;[b_-,b_+]}_{N,\beta}$
satisfies a large deviation principle with speed $N$ and rate function
equal to $\beta\widetilde{\mathcal{J}}^{V;[b_-,b_+]}_{\rm max}$
(resp. $\beta\widetilde{\mathcal{J}}^{V;[b_-,b_+]}_{\rm min}$)
on $[b_-,b_+]$. In other words, for any closed subset $F$, or open subset $\Omega$, of $[b_-,b_+]$:
\bea
\limsup_{N\ra\infty}\frac{1}{N}\ln \mu^{V;[b_-,b_+]}_{N,\beta}\left(\lambda_{\mathrm{max}}\in F\right) & \le & -\beta\,\inf_{x \in F} \widetilde{\mathcal{J}}^{V;[b_-,b_+]}_{\rm max}(x), \nn \\
\liminf_{N\ra\infty}\frac{1}{N}\ln \mu^{V;[b_-,b_+]}_{N,\beta}\left(\lambda_{\mathrm{max}} \in \Omega\right) & \ge & -\beta\,\inf_{x \in \Omega} \widetilde{\mathcal{J}}^{V;[b_-,b_+]}_{\rm max}(x), \nn
\eea
and similar statements hold for $\lambda_{\mathrm{min}}$.

In particular, for any $\varepsilon>0$,
\bea\label{bnmin}\limsup_{N\ra\infty}\frac{1}{N}\ln \mu^{V;[b_-,b_+]}_{N,\beta}\left(\lambda_{\mathrm{min}}\le \alpha_--\varepsilon \right) & <&0,\\
\label{bnmax}\limsup_{N\ra\infty}\frac{1}{N}\ln \mu^{V;[b_-,b_+]}_{N,\beta}\left(\lambda_{\mathrm{max}}\ge \alpha_++\varepsilon \right) & <&0.\eea
\end{itemize}
\end{proposition}
We give a proof of this proposition for completeness in the Appendix.
\subsection{Weak dependence on the soft edges}
We first state the global version of the result:
\begin{proposition}
\label{uuu2} Let $V\::\:[b_-,b_+] \rightarrow \mathbb{R}$ be a continuous function, and if $b_{\tau} = \tau \infty$, assume \eqref{assss}.
Suppose $b_- < \alpha_-$, and assume furthermore that the minimum value of $\mathcal{J}^{V;[b_-,b_+]}$ is achieved only on $[\alpha_-,\alpha_+]$.
Then, for any $\varepsilon > 0$, there exists $\eta_{\varepsilon} > 0$ so that:
\beq
Z_{N,\beta}^{V; [b_-,b_+]}=
Z^{V;[\alpha_- - \varepsilon,b_+]}_{N,\beta}\big(1+O(e^{-N\,\eta_\varepsilon})\big), \nn
\eeq
and there exists a universal constant $\gamma_n > 0$ such that, for any $x_1,\ldots,x_n\in (\mathbb{C}\backslash [b_-,b_+])^n$:
\beq
\label{eq:coinW} \big|W_n^{V;[b_-,b_+]}(x_1,\ldots,x_n) - W_n^{V;[\alpha_- - \varepsilon,b_+]}(x_1,\ldots,x_n)\big| \leq \frac{\gamma_n\,e^{-N\eta_{\varepsilon}}}{\prod_{i = 1}^{n} d(x_i,[b_-,b_+])}.
\eeq
A similar result holds for the upper edge.
 \end{proposition}
We also have a local version:

\begin{proposition}
\label{pro} Let $V\::\:[b_-,b_+] \rightarrow \mathbb{R}$ be a continuous function, and if $b_{\tau} = \tau \infty$, assume \eqref{assss}.
Suppose $b_- < \alpha_+$, and assume furthermore that the minimum value of $\mathcal{J}^{V;[b_-,b_+]}$ is achieved only on $[\alpha_-,\alpha_+]$. For any $\varepsilon > 0$ small enough, there exists $\eta_{\varepsilon} > 0$ so that, for any $a_- \in ]b_-,\alpha_- - \varepsilon[$:
\beq
\left|\partial_{a_-} \ln Z_{N}^{V;[a_-,b_+]}\right| \leq e^{-N\eta_{\varepsilon}}, \nn
\eeq
and, for any $x_1,\ldots,x_n \in (\mathbb{C}\setminus[a_-,b_+])^n$:
\beq
\forall x_1,\ldots,x_n \in \mathbb{C}\setminus [a_-,b_+],\quad
\left| \partial_{a_-} W_{n}^{V;[a_-,b_+]}(x_1,\ldots,x_n)\right| \leq \frac{\gamma_n\,N^{n}}{\prod_{i = 1}^{n} d(x_i,[a_-,b_+])}\,e^{-N\eta_{\varepsilon}}. \nn
\eeq
A similar statement holds for derivatives with respect to the upper bound.
\end{proposition}

\noindent {\bf Proof.}
If $b_- \neq \alpha_-$, let $a_-\in ]b_-,\alpha_-[$. Notice that:
\beq
\label{qho}\left(1 - \frac{Z_{N,\beta}^{V;[a_-,b_+]}}{Z_{N,\beta}^{V;[b_-,b_+]}}\right) = \mu_{N,\beta}^{V;[b_-,b_+]}[\lambda_{\mathrm{min}} \leq a_-].
\eeq
If now $\phi\::\: [b_-,b_+]^{N} \rightarrow \mathbb{C}$ is a bounded continuous function, we can write:
{\small \beq
\mu_{N,\beta}^{V;[b_-,b_+]}[\phi(\lambda)] - \mu_{N,\beta}^{V;[a_-,b_+]}[\phi(\lambda)] = \mu_{N,\beta}^{V;[b_-,b_+]}\big[\phi(\lambda)\,\mathbf{1}_{\lambda_{\mathrm{min}} \leq a_-}\big] + \left(\frac{Z_{N,\beta}^{V;[a_-,b_+]}}{Z_{N,\beta}^{V;[b_-,b_+]}} - 1\right)\mu_{N,\beta}^{V;[a_-,b_+]}[\phi(\lambda)]. \nn
\eeq}
$\!\!\!$Thus, we find:
\bea
\big|\mu_{N,\beta}^{V;[b_-,b_+]}[\phi(\lambda)] - \mu_{N,\beta}^{V;[a_-,b_+]}[\phi(\lambda)]\big|
& \leq & 2\,\big(\mathrm{sup}_{\lambda \in [b_-,b_+]^N} |\phi(\lambda)|\big)\,\mu_{N,\beta}^{V;[b_-,b_+]}[\lambda_{\mathrm{min}} \leq a_-]. \nn
\eea
This can be applied for the disconnected correlators:
\beq
\overline{W}_n^{V;[b_-,b_+]}(x_1,\ldots,x_n) =  \mu_{N,\beta}^{V;[b_-,b_+]}\Big[\prod_{j = 1}^{n} \sum_{i_j = 1}^{N} \frac{1}{x_j - \lambda_{i_j}} \Big], \nn
\eeq
and we obtain:
\beq
\big|\overline{W}_n^{V;[b_-,b_+]}(x_1,\ldots,x_n) - \overline{W}_n^{V;[a_-,b_+]}(x_1,\ldots,x_n)\big| \leq \frac{2\,N^n}{\prod_{j = 1}^{n} d(x_j,[b_-,b_+])}\,\mu_{N,\beta}^{V;[b_-,b_+]}[\lambda_{\mathrm{min}} \leq a_-]. \nn
\eeq
Similarly, one finds:
\beq
\big|\partial_{a_-} W_n^{V;[a_-,b_+]}(x_1,\ldots,x_n) \big| \leq \frac{2\,N^{n}}{\prod_{j = 1}^{n} d(x_j,[a_-,b_+])}\,\partial_{a_-}\ln Z_{N,\beta}^{V;[a_-,b_+]}. \nn
\eeq
The correlators $W_n^{V;[b_-,b_+]}$ are just sums of monomials of the form $W_{n_1}^{V;[b_-,b_+]}(I_1)\cdots W_{n_m}^{V;[b_-,b_+]}(I_m)$ where $I_1,\ldots,I_m$ is a partition of $\{x_1,\ldots,x_n\}$.
So, it is enough to establish the weak dependence at the level of the partition function. The global version is a direct consequence of Eqn.~\ref{bnmin} applied to Eqn.~\ref{qho}:
\beq
\limsup_{N \rightarrow \infty} \frac{1}{N}\ln\left(1 - \frac{Z_{N,\beta}^{V;[a_-,b_+]}}{Z_{N,\beta}^{V;[b_-,b_+]}}\right)  < 0. \nn
\eeq
For the local version, we rather need to bound:
\beq
\partial_{a_-} \ln Z_{N,\beta}^{V;[a_-,b_+]} = N\,\frac{Z_{N - 1,\beta}^{\frac{NV}{N - 1};[a_-,b_+]}}{Z_{N,\beta}^{V;[a_-,b_+]}}\,\mu_{N - 1,\beta}^{\frac{NV}{N - 1};[a_-,b_+]}\left[e^{\beta\left(-\frac{NV(a_-)}{2} + \sum_{i = 1}^{N - 1} \ln|\lambda_i - a_-|\right)}\right]. \nn
\eeq
If $a_- \in ]b_-,\alpha_-[$ is fixed, by the large deviation principle for $L_{N -1}$ under $\mu_{N - 1}^{\frac{NV}{N - 1};[a_-,b_+]}$, since the logarithm is a lower semicontinuous function,
there exists $\gamma > 0$ such that, for any $\epsilon > 0$, for $N$ large enough:
\beq
\mu_{N - 1,\beta}^{\frac{NV}{N - 1};[a_-,b_+]}\left[e^{\beta\left(-\frac{NV(a_-)}{2} +\sum_{j = 1}^{N - 1} \ln|a_- - \lambda_j|\right)}\right] \leq \gamma\,e^{-\beta N(1 - \epsilon) {\mathcal{J}}^{V;[b_-,b_+]}(a_-)}. \nn
\eeq
Moreover, we have seen in Eqn.~\ref{cont3} that for $N$ large enough:
\beq
\frac{Z_{N - 1,\beta}^{\frac{NV}{N - 1};[a_-,b_+]}}{Z_{N,\beta}^{V;[a_-,b_+]}}\le
e^{\beta N(1 - \epsilon)\inf_{\xi\in [b_-,\alpha_-]} {\mathcal{J}}^{V;[a_-,b_+]}_{\rm min}(\xi)}. \nn
\eeq
By assumption, $\widetilde{\mathcal{J}}_{\mathrm{min}}^{V;[b_-,b_+]}(a_-) = \mathcal{J}_{\mathrm{min}}^{V;[b_-,b_+]}(a_-) - \inf_{\xi \in [b_-,\alpha_-]} \mathcal{J}^{V;[b_-,b_+]}(\xi) > 0$, leading to:
\beq
\limsup_{N \rightarrow \infty} \frac{1}{N}\,\ln \left|\partial_{a_-} \ln Z_{N,\beta}^{V;[a_-,b_+]}\right| < 0, \nn
\eeq
which is the bound we sought. The arguments at the upper edge are similar. \hfill $\diamond$

\section{Loop equations}\label{secloop}
We shall assume in this Section and also in Section~\ref{sec:rec1} that the $\lambda_i$ are integrated over a segment $[a_-,a_+]$ with: 
\begin{hypothesis}
\label{hypq}
$-\infty < a_- < a_+ < +\infty$.
\end{hypothesis}
Indeed, considering finite intervals $[a_-,a_+]$ is convenient to ensure from the beginning that the Cauchy-Stieltjes transform
yields functions which are holomorphic outside $[a_-,a_+]$. We also assume in this section:
\begin{hypothesis}
\label{Vhyp}
$V\,:\, [a_-,a_+] \rightarrow \mathbb{C}$ can be extended as a holomorphic function in some open neighborhood of $[a_-,a_+]$.
\end{hypothesis}
This will allow us to use complex analysis (Cauchy residue formula, moving the contours, etc.)

We shall derive  the "loop equations", also called Schwinger-Dyson equations or Pastur equations \cite{Pastur72} in this context. These equations express the invariance by change of variable of an integration, up to boundary terms. We stress that these equations are exact for finite $N$. Although the technique is well-known, we recall the derivation here for the $\beta$ matrix models with edges $a_-,a_+$ in order to have a self-contained presentation.

\subsection{First version}

\begin{theorem}
\label{th:loopeqt} Loop equation at rank $1$. For any $x \in \mathbb{C}\setminus[a_-,a_+]$:
\bea
W_{2}(x,x) + \big(W_1(x)\big)^2 + \Big(1 - \frac{2}{\beta}\Big)\frac{\mathrm{d}}{\mathrm{d}x}\Big(W_1(x)\Big)+ \frac{N(1 - 2/\beta) - N^2}{(x - a_-)(x - a_+)} & & \nn \\
 - N\left(\oint_{\mathcal{C}([a_-,a_+])}\frac{\mathrm{d}\xi}{2i\pi}\,\frac{1}{x - \xi}\,\frac{(\xi - a_-)(\xi - a_+)}{(x - a_-)(x - a_+)}\,V'(\xi)\,W_1(\xi)\right) & = & 0. \nn
\eea
$\mathcal{C}([a_-,a_+])$ is a contour surrounding $[a_-,a_+]$ in positive orientation, and included in the domain where $V'$ is holomorphic.
\end{theorem}

\begin{theorem}
\label{th:loopeqtn} Loop equation at rank $n$. Let $x_I = (x_i)_{i \in I}$ be a $(n - 1)$-uple of spectator variables in $(\mathbb{C}\setminus[a_-,a_+])^{n - 1}$. For any $x \in \mathbb{C}\setminus [a_-,a_+]$:
\bea
W_{n + 1}(x,x,x_I) + \sum_{J \subseteq I} W_{|J| + 1}(x,x_J)\,W_{n - |J|}(x,x_{I\setminus J}) + \Big(1 - \frac{2}{\beta}\Big)\frac{\mathrm{d}}{\mathrm{d} x}\Big(W_n(x,x_I)\Big) && \nonumber \\
 - N\left(\oint_{\mathcal{C}([a_-,a_+])} \frac{\mathrm{d}\xi}{2i\pi}\,\frac{1}{x - \xi}\,\frac{(\xi - a_-)(\xi - a_+)}{(x - a_-)(x - a_+)}\,V'(\xi)\,W_n(\xi,x_I)\right) & & \nonumber \\
+ \frac{2}{\beta}\sum_{i \in I} \frac{\mathrm{d}}{\mathrm{d} x_i}\Big(\frac{W_{n - 1}(x,x_{I\setminus\{i\}}) - \frac{(x_i - a_-)(x_i - a_+)}{(x - a_-)(x - a_+)}\,W_{n - 1}(x_I)}{x - x_i}\Big) & = & 0. \nonumber
\eea
\end{theorem}

\noindent {\bf  Proof of Theorem \ref{th:loopeqt}.}
For any smooth real-valued function $h$,  and $\varepsilon > 0$ small enough,
\beq
\psi_{h,\varepsilon}\,:\, \lambda \mapsto \lambda + \varepsilon h(\lambda) \nn
\eeq
defines a differentiable family of diffeomorphisms from $[a_-,a_+]$ to some interval $\psi_{h,\varepsilon}([a_-,a_+])$.
We assume hereafter that $h(a_-)=h(a_+)=0$ so that $\psi_{h,\varepsilon}([a_-,a_+])=[a_-,a_+]$ for $\varepsilon$ small enough. We have:
\beq
\label{eq:inv}1= \int_{[a_-,a_+]^N} \mathrm{d}\mu^V_{N,\beta}\big(\psi_{h,\varepsilon}(\lambda_1),\ldots,\psi_{h,\varepsilon}(\lambda_N)\big).
\eeq
When $\varepsilon \rightarrow 0$, the first subleading order of the right hand side must vanish. It can be computed in three parts. A first term comes from the variation of the Lebesgue measure $\prod_{i} \mathrm{d}\lambda_i$, which is given by the Jacobian of the change of variable:
\beq
\Big(\prod_{i = 1}^N \mathrm{d}\psi_{h,\varepsilon}(\lambda_i)\Big) =
\Big(\prod_{i = 1}^N \mathrm{d}\lambda_i\Big)\Big(1 + \varepsilon
\int h'(\xi)\dd M_N(\xi) + o(\varepsilon)\Big). \nn
\eeq
A second term comes from the variation of the Vandermonde:
\bea
&&|\Delta\big(\psi_{h,\varepsilon}(\lambda)\big)|^{\beta} =  |\Delta(\lambda)|^{\beta}\Big[1 + \varepsilon\,\beta\sum_{1 \leq i < j \leq N} \frac{h(\lambda_i) - h(\lambda_j)}{\lambda_i - \lambda_j} + o(\varepsilon)\Big] \nonumber \\
& = & |\Delta(\lambda)|^{\beta}\left\{1 +
\varepsilon\,\frac{\beta}{2}\Big(
\iint \frac{h(\xi)-h(\eta)}{\xi-\eta} \dd M_N(\xi)\dd M_N(\eta) -\int h'(\xi)\dd M_N(\xi)\Big) + o(\varepsilon)\right\}. \nonumber
\eea
The last term comes from the variation of the Boltzmann weight:
\beq
\prod_{i = 1}^N
e^{-\frac{N\beta}{2}\,V[\psi_{h,\varepsilon}(\lambda_i)]} =
\Big(\prod_{i = 1}^N e^{-\frac{N\beta}{2}\,V(\lambda_i)}\Big)\Big(1 - \varepsilon\,\frac{N\beta}{2} \int V'(\xi)\,h(\xi)\dd M_N(\xi) + o(\varepsilon)\Big). \nn
\eeq
Summing all terms up, the first order in $\varepsilon$ in Eqn.~\ref{eq:inv} vanishes iff:
\bea
 &&\mu_{N,\beta}^{V;[a_-,a_+]}\Big[\iint \frac{h(\xi) - h(\eta)}{\xi - \eta}\,\dd M_N(\xi)\dd M_N(\eta)
 - N \int V'(\xi)h(\xi)\,\dd M_N(\xi)\Big] \nn \\
&=&
\label{bbb}\Big(1-\frac{2}{\beta} \Big)\mu_{N,\beta}^{V;[a_-,a_+]}\Big[\int h'(\xi)\,\dd M_N(\xi) \Big].
\eea
Note that even though this equation was obtained for real-valued
functions $h$, we can at this point remove this condition by linearity. To obtain an equation involving correlators, one can take for $x \in \mathbb{C}\setminus[a_-,a_+]$ the function $h$ defined by:
\beq
h(\xi) = \frac{(\xi - a_-)(\xi - a_+)}{x - \xi} = \frac{(x - a_-)(x - a_+)}{x - \xi} + a_- + a_+ - x - \xi.\nn
\eeq
thus preserving $[a_-,a_+]$. We recall that $V'$ is holomorphic in a neighborhood of $[a_-,a_+]$. So, by Cauchy formula,
for any contour $\mathcal{C}([a_-,a_+])$ surrounding $[a_-,a_+]$ inside this neighborhood and not enclosing $x$:
\beq
\int \frac{V'(\xi)(\xi-a_-)(\xi-a_+)}{x-\xi}\,\dd M_N(\xi) =
\int \dd M_N(\xi)\oint_{\mathcal{C}([a_-,a_+])}\frac{\dd\eta}{2i\pi}\,\frac{V'(\eta)(\eta-a_-)(\eta-a_+)}{(\eta - \xi)(x - \eta)}. \nn
\eeq
Hence, we obtain:
\bea
\label{eq:loopB} && W_2(x,x) + \big(W_1(x)\big)^2 - \frac{N^2}{(x - a_-)(x - a_+)} \nn \\
&& - N\oint_{\mathcal{C}([a_-,a_+])} \frac{\dd\xi}{2i\pi}\,\frac{1}{x - \xi}\,\frac{(\xi - a_-)(\xi - a_+)}{(x - a_-)(x - a_+)}\,V'(\xi)\,W_1(\xi)  \nn \\
& = & \Big(1 - \frac{2}{\beta}\Big)\Big(-\frac{\dd}{\dd x}\big(W_1(x)\big) - \frac{N}{(x - a_-)(x - a_+)}\Big). \nn
\eea

\noindent
{\bf  Proof of  Theorem \ref{th:loopeqtn}.}
By definition of the cumulants, if we define a shifted potential $V_{(x;\epsilon)}(\xi) = V(\xi) + \frac{\epsilon}{x - \xi}$, we have:
\beq
W_{n}^{V}(x,x_2,\ldots,x_n) = -\frac{2}{\beta N} \partial_{\epsilon}\Big(W_{n - 1}^{V_{(x;\epsilon)}}(x_2,\ldots,x_{n})\Big)\Big|_{\epsilon = 0} .\nn
\eeq
Notice that the matrix integral with this shifted potential is still convergent, because the eigenvalues live on the finite interval $[a_-,a_+]$.
Therefore, we can obtain the loop equations at rank $n$ by taking a perturbed potential in Thm.~\ref{th:loopeqt}:
\beq
V_{(x_2;\epsilon_2),\ldots,(x_{n};\epsilon_{n})}(\xi) = V(\xi) + \sum_{i = 2}^n \frac{\epsilon_i}{x_i - \xi}, \nn
\eeq
and identifying the term in $\Big[\prod_{i = 2}^{n} \big(\frac{-2}{\beta N}\big) \epsilon_i\Big]$ when $\epsilon_i \rightarrow 0$.

\subsection{Second version}

Here is another equivalent form of the loop equations. All $W_n$ depend implicitly on the interval of integration $[a_-,a_+]$.

\begin{theorem}
\label{th:loopeqtt} Loop equation at rank $1$. For any $x \in \mathbb{C}\setminus[a_-,a_+]$:
\bea
W_{2}(x,x) + \big(W_1(x)\big)^2 + \Big(1 - \frac{2}{\beta}\Big)\frac{\mathrm{d}}{\mathrm{d}x}\Big(W_1(x)\Big) & & \nn \\
- N \left(\oint_{\mathcal{C}([a_-,a_+])} \frac{\dd\xi}{2i\pi}\,\frac{V'(\xi)\,W_1(\xi)}{x - \xi}\right) - \frac{2}{\beta}\left(\frac{\partial_{a_-}\ln Z}{x - a_-} + \frac{\partial_{a_+} \ln Z}{x - a_+}\right) & = & 0. \nn
\eea
$\mathcal{C}([a_-,a_+])$ is a contour surrounding $[a_-,a_+]$ in positive orientation, and included in the domain where $V$ is holomorphic.
\end{theorem}

\begin{theorem}
\label{th:loopeqttn} Loop equation at rank $n$. Let $x_I = (x_i)_{i \in I}$ a $(n - 1)$-uple of spectator variables in $(\mathbb{C}\setminus[a_-,a_+])^{n - 1}$. For any $x \in \mathbb{C}\setminus\mathbb[a_-,a_+]$:
\bea
W_{n + 1}(x,x,x_I) + \sum_{J \subseteq I} W_{|J| + 1}(x,x_J)\,W_{n - |J|}(x,x_{I\setminus J}) + \left(1 - \frac{2}{\beta}\right)\frac{\mathrm{d}}{\mathrm{d} x}\Big(W_n(x,x_I)\Big) && \nonumber \\
 - N\left(\oint_{\mathcal{C}([a_-,a_+])} \frac{\mathrm{d}\xi}{2i\pi}\,\frac{V'(\xi)\,W_n(\xi,x_I)}{x - \xi}\right) & & \nonumber \\
+ \frac{2}{\beta}\sum_{i \in I} \frac{\mathrm{d}}{\mathrm{d} x_i}\left(\frac{W_{n - 1}(x,x_{I\setminus\{i\}}) - W_{n - 1}(x_I)}{x - x_i}\right)  - \frac{2}{\beta}\left(\frac{\partial_{a_-} W_{n - 1}(x_I)}{x - a_-} + \frac{\partial_{a_+} W_{n - 1}(x_I)}{x - a_{+}}\right) & = & 0. \nonumber
\eea
\end{theorem}

\noindent \textbf{Proof}
In the former proof, if we use a change of variable $h$ which does not preserve $[a_-,a_+]$, the partition function becomes (to first order in $\varepsilon$):
{\small \beq
Z_{N}^{V;\psi_{h,\varepsilon}([a_-,a_+])} \rightarrow Z_{N}^{V;[a_-,a_+]}\left[1 + \varepsilon\left(h(a_-)\,\partial_{a_-} \ln Z_{N}^{V;[a_-,a_+]} + h(a_+)\,\partial_{a_+} \ln Z_{N}^{V;[a_-,a_+]}\right) + o(\varepsilon)\right]. \nn
\eeq}
$\!\!\!$Thus, Eqn.~\ref{bbb} receives those extra terms, and becomes:
{\small\bea
 &&\mu_{N,\beta}^{V;[a_-,a_+]}\Big[\iint \frac{h(\xi) - h(\eta)}{\xi - \eta}\,\dd M_N(\xi)\dd M_N(\eta)
 - N \int V'(\xi)h(\xi)\,\dd M_N(\xi)\Big] \nn \\
 &=&
\Big(1-\frac{2}{\beta} \Big)\mu_{N,\beta}^{V;[a_-,a_+]}\Big[\int h'(\xi)\,\dd M_N(\xi) \Big] +  \frac{2}{\beta}\left(h(a_-)\,\partial_{a_-}\ln Z_{N}^{V;[a_-,a_+]} + h(a_+)\,\partial_b \ln Z_{N}^{V;[a_-,a_+]}\right). \nn
\eea}
$\!\!\!$In particular, when we choose $h(\xi) = \frac{1}{x - \xi}$, we obtain:
\bea
& & W_2(x,x) + \big(W_1(x)\big)^2 - N\oint_{\mathcal{C}([a_-,a_+])} \frac{\dd\xi}{2i\pi}\,\frac{V'(\xi)W_1(\xi)}{x - \xi} \nn \\
& = & -\left(1 - \frac{2}{\beta}\right)\frac{\dd}{\dd x}\Big(W_1(x)\Big) + \frac{2}{\beta}\,\frac{\partial_{a_-} \ln Z}{x - a_-} + \frac{2}{\beta}\,\frac{\partial_{a_+} \ln Z}{x - a_+}. \nn
\eea
The loop equation at higher rank can be deduced as before by perturbing the potential. \hfill$\diamond$

\subsection{Remark}
\label{sechy}If we compare those expressions to the first version of the loop equations, we find by consistency:
\beq
\partial_{a_\tau} \ln Z = \frac{1}{a_{-\tau} - a_{\tau}}\,\left\{-\frac{N^2\beta}{2} + N\left(\frac{\beta}{2} - 1\right) + \frac{N\beta}{2} \oint_{\mathcal{C}([a_-,a_+])} \frac{\dd\xi}{2i\pi}\,(\xi - a_{-\tau})V'(\xi)W_1(\xi)\right\}, \nn
\eeq
and for higher correlators $\partial_{a_\tau} W_{n - 1}(x_I)$ equals
\beq
  \frac{1}{a_{-\tau} - a_{\tau}}\,\left\{\frac{N\beta}{2}\oint_{\mathcal{C}([a_-,a_+])} \frac{\dd\xi}{2i\pi}\,(\xi - a_{-\tau})V'(\xi) W_n(\xi) + \sum_{i \in I} \frac{\dd}{\dd x_i}\Big((x_i - a_{-\tau})W_{n - 1}(x_I)\Big)\right\} \nn
\eeq
for $\tau \in \{\pm\}$.

\section{The $1/N$ expansion}
\label{sec:rec1}
\subsection{Notations, assumptions, proposition}

This section relies on complex analysis and inequalities for probability measures. We make four assumptions on the potential $V$, which are valid only in this section.
The link with our main theorem will be done in Section~\ref{final}.

We keep on with the assumption:
\label{hypotheses}
\begin{hypothesis}
\label{hypq2}
$-\infty < a_- < a_+ < +\infty$.
\end{hypothesis}

\vspace{0.2cm}

Since $V$ is smooth, the equilibrium measure $\dd \mu_{\mathrm{eq}}^{V;[a_-,a_+]}(\xi)$ will in fact be a density $\rho(\xi)\dd\xi$, where $\rho\,:\, [a_-,a_+] \rightarrow [0,+\infty]$ is a continuous function. We call $\mathrm{supp}\,\rho = \overline{\{x \in [a_-,a_+]\quad \rho(x) > 0\}}$ its support. In the hermitian case ($\beta = 2)$, a $1/N$ expansion is expected only when $\mathrm{supp}\,\rho$ is connected. We assume here also:
\begin{hypothesis}
\label{onecut}$V$ leads to a one-cut regime, i.e. the support of $\mu_{\mathrm{eq}}^{V;[a_-,a_+]}$ is an interval $[\alpha_-,\alpha_+] \subseteq [a_-,a_+]$.
\end{hypothesis}

\vspace{0.2cm}

In order to write the loop equations as in Section~\ref{secloop}, we assume:
\begin{hypothesis}
\label{hypVreal}
$V$ is real-valued on $[a_-,a_+]$, and can be extended as a holomorphic function on some open neighborhood $U$ of $[a_-,a_+]$.
\end{hypothesis}

\vspace{0.2cm}
We justify in Remark~\ref{Rem} later that there exists a unique analytic function $y\,:\,U \rightarrow \mathbb{C}\,\cup\,\{\infty\}$ such that, for any $x \in [\alpha_-,\alpha_+]$, we have:
\beq
\label{rhom}\rho(x) = \frac{1}{i\pi}\,\lim_{\epsilon \rightarrow 0^+} y(x + i\epsilon).
\eeq
This function can be written $y(x) = S(x)\sigma(x)$, where $S$ is now a holomorphic function defined on $U$, and $\sigma$ is of the form:
\beq
\label{sigm}\sigma(x) = \sqrt{\frac{\prod_{\tau \in \mathsf{Soft}} (x - \alpha_\tau)}{\prod_{\tau' \in \mathsf{Hard}} (x - \alpha_{\tau'})}}.
\eeq
The lower edge $a_-$ is
\begin{itemize}
\item[$\bullet$] either a hard edge, meaning that $a_- = \alpha_-$. Then, $\rho(x) \in O\big((x - \alpha_-)^{-1/2}\big)$ when $x \rightarrow \alpha_-$.
\item[$\bullet$] or a soft edge, meaning that $a_- < \alpha_-$. Then, $\rho(x) \in O\big((x - \alpha_-)^{1/2}\big)$ when $x \rightarrow \alpha_-$.
\end{itemize}
and the same distinction exists independently for the upper edge $a_+$. Our discussion holds for both hard and soft cases. However, a key technical assumption is:
\begin{hypothesis}
\label{offc} $V$ is offcritical on $[a_-,a_+]$, in the sense that $S(x)$ remains positive on $[a_-,a_+]$.
\end{hypothesis}
For instance, Hyp.~\ref{onecut} and \ref{offc} automatically hold when $V$ is strictly convex. For a generic $V$ satisfying Hyp.~\ref{onecut}, we have $S(\alpha_-) > 0$ and $S(\alpha_+) > 0$, so we can always find an interval $[a_-,a_+]$ which is a strict enlargment of $[\alpha_-,\alpha_+]$, such that Hyp.~\ref{offc} holds on $[a_-,a_+]$. We call "critical point on $[a_-,a_+]$", the situation corresponding to a choice of $V$ such that $S$ has a zero on $[a_-,a_+]$. In this article, we do not tackle the question of the double scaling limit for $\beta$ matrix models ($N \rightarrow +\infty$ and coefficients of $V$ finely tuned with $N$ to achieve a critical point when $N = \infty$). Though, this would be a very interesting regime in relation with universality questions, considering the absence of Riemann-Hilbert techniques when $\beta \neq 1,2,4$.

\begin{figure}[h]
\begin{center}\includegraphics[width = 0.8\textwidth]{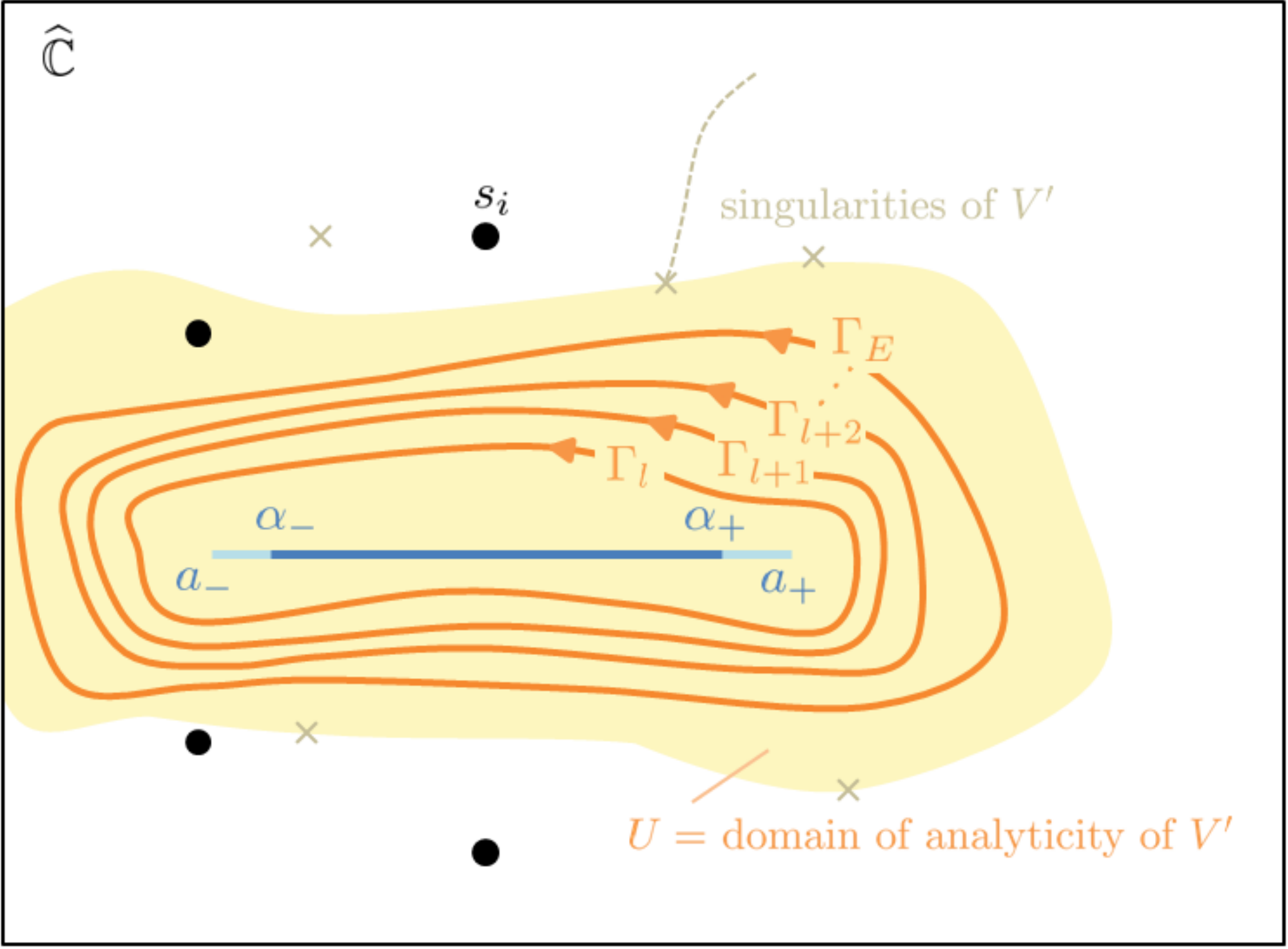}
\caption{\label{fig:cont1} Hypothesis on the location of $s_i$ and contour configurations}
\end{center}
\end{figure}

We shall allow $V$ itself to depend on $N$ and have a $1/N$ expansion. To give precise statements about those expansions, we need some notations. For any Jordan curve $\Gamma$, we note $\mathrm{Ext}(\Gamma)$ (resp. $\mathrm{Int}(\Gamma)$) the unbounded (resp. bounded) connected component of $\mathbb{C}\setminus\Gamma$. In the following, we fix once for all a Jordan curve $\Gamma_E$, and a sequence of nested Jordan curves $(\Gamma_l)_{l \in \mathbb{N}}$, which all live in $\mathbb{C}\setminus[a_-,a_+]$, and such that
\begin{itemize}
\item[$(i)$] $\Gamma_E \subseteq U$.
\item[$(ii)$] $\{x \in U\quad S(x) = 0\}\cap \mathrm{Int}(\Gamma_E) = \emptyset$.
\item[$(iii)$] $\forall l \in \mathbb{N}\quad  \Gamma_l \subseteq \mathrm{Int}(\Gamma_{l + 1})$.
\item[$(iv)$] $\forall l \in \mathbb{N} \quad \Gamma_l \subseteq \mathrm{Int}(\Gamma_E)$.
\end{itemize}
The contour configuration is depicted in Fig.~\ref{fig:cont1}, where the zeroes of $S$ were called $s_i$. In the remaining of the text, $\Gamma$ will refer to a Jordan curve in $\mathrm{Int}(\Gamma_E)\setminus[a_-,a_+]$. We will use the following norm on the space of continuous functions on a contour $\Gamma$
$$\p f\p_\Gamma= \sup_{x_i \in \Gamma} |f(x_1,\ldots,x_n)|.$$
On the space $\mathcal{H}_{n;[a_-,a_+]}^{(1)}$ of holomorphic functions on $(\mathbb{C}\setminus[a_-,a_+])^n$, which behave as $O(1/x_i)$ when $x_i \rightarrow \infty$, we have by the maximum principle
\beq
\p f \p_{\Gamma}= \sup_{x_i \in \mathrm{Ext}(\Gamma)} |f(x_1,\ldots,x_n)|. \nn
\eeq
 One can easily derive the following useful inequalities:
\beq
\label{eq:contr}\forall f \in \mathcal{H}_{1;[a_-,a_+]}^{(1)}\quad \forall x_0\in \mathrm{Ext}(\Gamma_{l+1}) 
\quad  \forall l \in \mathbb{N} \qquad \Big\| \frac{f(\bullet) - f(x_0)}{\bullet - x_0} \Big\|_{\Gamma_l} \leq \p f' \p_{\Gamma_{l + 1}} \leq \zeta_{l}\,\p f \p_{\Gamma_l},
\eeq
where $\zeta_{l} = \frac{\ell(\Gamma_l)}{2\pi\,d^2(\Gamma_l,\Gamma_{l + 1})}$ is a finite constant depending only on the relative position of $\Gamma_l$ and $\Gamma_{l  + 1}$.

\vspace{0.1cm}

\vspace{0.1cm}

\noindent Now, we can state our last assumption:
\begin{hypothesis}
\label{hypVN}
$V$ admits a $1/N$ asymptotic expansion:
\beq
V(x) \,\mathop{=}_{\Gamma_E}\, \sum_{k \geq 0} N^{-k}\,V^{\{k\}}(x), \nn
\eeq
with functions $V^{\{k\}}$ independent of $N$, such that $V^{\{k\}}$ is holomorphic in $U$. The $\mathop{=}_{\Gamma_E}$ equality means that, for any positive integer $K$, there exists a positive constant $v_K$ such that, for $N$ large enough:
\beq
\mathop{\mathrm{sup}}_{\xi \in \Gamma_E} \Big|V(\xi) - \sum_{k = 0}^{K} N^{-k}\,V^{\{k\}}(\xi)\Big| \leq N^{-(K + 1)}\,v_K .\nn
\eeq
(The maximum principle implies automatically the same bound with $\Gamma$ replacing $\Gamma_E$ as $V$ is analytic in $\mathrm{Int}(\Gamma_E)$). 
\end{hypothesis}
In many applications, $V$ is independent of $N$ (i.e. $V \equiv
V^{\{0\}}$). There is however no difficulty in our reasoning to consider
potentials which depend on $N$ within Hyp.~\ref{hypVN}.

\vspace{0.1cm}

Our intermediate result is:

\begin{theorem}
\label{recerror}
If Hyp.~\ref{hypq2}-\ref{hypVN} hold, the correlators admit an asymptotic expansion when $N \rightarrow \infty$ with respect to the norm $\p \cdot \p_{\Gamma_E}$, of the form:
\beq
\forall n \geq 1,\qquad W_n = \sum_{k \geq n - 2} N^{-k}\,W_n^{\{k\}}, \nn
\eeq
where $W_n^{\{k\}} \in \mathcal{H}_{n;[\alpha_-,\alpha_+]}^{(1)}$.
\end{theorem}

\subsection{Relevant linear operators}

\subsubsection{The operator $\mathcal{K}$}

\label{opK}
We introduce the following linear operator defined on the space $\mathcal{H}^{(2)}_{1;[a_-,a_+]}$ of holomorphic functions on $\mathbb{C}\setminus[a_-,a_+]$ which behave as $O(1/x^2)$ when $x \rightarrow \infty$:
\beq
(\mathcal{K}f)(x) = 2W_1^{\{-1\}}(x)f(x) - \oint_{\mathcal{C}([a_-,a_+])}\frac{\mathrm{d}\xi}{2i\pi}\,\frac{L(\xi)}{L(x)}\Big(\frac{1}{x - \xi} + c\Big)\,(V^{\{0\}})'(\xi)\,f(\xi). \nn
\eeq
This operator for an appropriate choice of $L$ and $c$ appears in the loop equations. We have found the following choice convenient:
\bea
L(x) & = & \prod_{\tau \in \mathsf{Hard}} (x - a_\tau).\quad \nn \\
c & = & \left\{\begin{array}{ll} 0 & \mathrm{if}\,\,\mathsf{Soft} = \{\pm\}\,\,\mathrm{or}\,\,\mathsf{Hard} = \{\pm\}, \\ \frac{1}{a_{\tau} - a_{-\tau}} & \mathrm{if}\,\,\tau \in \mathsf{Soft}\,\,\mathrm{and}\,\,(-\tau) \in \mathsf{Hard}.\end{array}\right. \nn
\eea
We may also rewrite:
\beq
\label{formuu} (\mathcal{K}f)(x) =  -2y(x)\,f(x) + \frac{(\mathcal{Q}\,f)(x)}{L(x)},
\eeq
with:
\beq
(\mathcal{Q}\,f)(x) = -\oint_{\mathcal{C}([a_-,a_+])\cup\mathcal{C}(x)} \frac{\dd \xi}{2i\pi}\,L(\xi)\Big(\frac{1}{x - \xi} + c\Big)(V^{\{0\}})'(\xi)\,f(\xi), \nn
\eeq
where $\mathcal{C}(x)$ is a contour surrounding $x$ only (computing a residue at $x$). For any $f \in \mathcal{H}_{1;[a_-,a_+]}^{(1)}$, $(\mathcal{Q}\,f)$ is analytic, with singularities only where $(V^{\{0\}})'$ has singularities, in particular is holomorphic in the neighborhood of $[a_-,a_+]$. We have set:
\beq
y(x) = -W_1^{\{-1\}}(x) + \frac{(V^{\{0\}})'(x)}{2}. \nn
\eeq
$y$ is discontinuous on the support of $\mu_{\mathrm{eq}}^{V;[a_-,a_+]}$ (see Thm.~\ref{th:1}), i.e. on $[\alpha_-,\alpha_+] \subseteq [a_-,a_+]$, but analytic on $\mathbb C\setminus[\alpha_-,\alpha_+]$. We justify in Remark~\ref{Rem} that $y(x) = S(x)\sigma(x)$ where $\sigma(x)$ was introduced in Eqn.~\ref{sigm} and the squareroot is chosen with its usual discontinuity on $\mathbb{R}_{-}$. Let us call $s_i \neq \alpha_-,\alpha_+$ the zeroes of $S(x)$ in the complex plane, and we assume that they do not lie in $[a_-,a_+]$ (Hyp.~\ref{offc}).

It is clear that $\mathrm{Im}\,\mathcal{K} \subseteq \mathcal{H}_{1;[a_-,a_+]}^{(1)}$. Here, $W_1^{\{-1\}}$ (hence $y$) has only cut $[\alpha_-,\alpha_+]$, and this operator is invertible\footnote{In general, on the space of holomorphic functions with $g + 1$ cuts, $\mathrm{dim}\,\mathrm{Ker}\,\mathcal{K} = g$, and one has to prescribe $g$ cycle integrals in order to define an inverse operator.}. Its inverse can be explicitly written, it is given by Tricomi formula \cite{Tricomi}:
\begin{equation}\label{eq:Km1}
 \forall x \in \mathbb{C}\setminus[a_-,a_+],\, \forall g \in \mathrm{Im}\,\mathcal{K},\quad \big(\mathcal{K}^{-1}g\big)(x) = \oint_{{\mathcal C}([a_-,a_+])} \frac{\mathrm{d}\xi}{2i\pi}\,\frac{1}{\xi - x}\,\frac{\widetilde{\sigma}(\xi)}{\widetilde{\sigma}(x)}\,\frac{g(\xi)}{2y(\xi)}. \end{equation}
where $\widetilde{\sigma}(x) = \sqrt{(x - \alpha_-)(x - \alpha_+)}$, and where we integrate over a contour surrounding $[a_-,a_+]$ but not $x$. Indeed, if $g \in \mathrm{Im}\,\mathcal{K}$, we can write for any $x \in \mathbb{C}\setminus [a_-,a_+]$:
\bea
\widetilde{\sigma}(x)\,f(x) & = & \Res_{\xi \rightarrow x} \frac{\mathrm{d}\xi}{\xi - x}\,\widetilde{\sigma}(\xi)\,f(\xi) \nonumber \\
& = & - \oint_{\mathcal{C}([a_-,a_+])} \frac{\mathrm{d}\xi}{2i\pi}\,\frac{\widetilde{\sigma}(\xi)\,f(\xi)}{\xi - x} \nonumber \\
& = & - \oint_{\mathcal{C}([a_-,a_+])} \frac{\mathrm{d}\xi}{2i\pi}\,\frac{\widetilde{\sigma}(\xi)}{\xi - x}\,\frac{1}{2y(\xi)}\left(-g(\xi) + \frac{(\mathcal{Q}f\big)(\xi)}{L(\xi)}\right) \nonumber \\
& = & \oint_{\mathcal{C}([a_-,a_+])} \frac{\mathrm{d}\xi}{2i\pi}\,\frac{\widetilde{\sigma}(\xi)}{\xi - x}\,\frac{g(\xi)}{2y(\xi)}. \nonumber
\eea
In the second line, we moved the contour from a neighborhood of $x$ to a neighborhood of $[a_-,a_+]$,
 and used the fact that $\widetilde{\sigma}(\xi) \in O(\xi)$ and $f(\xi) \in O(1/\xi^2)$, so that the residue at $\infty$ of the integrand vanishes.
In the fourth line, we use the fact that $L$ is chosen such that $\frac{\widetilde{\sigma}(\xi)}{y(\xi)L(\xi)} = \frac{1}{S(\xi)}$, which is holomorphic in a neighborhood of $[a_-,a_+]$ thanks to Hyp.~\ref{offc}. Since  $\big(\mathcal{Q}f\big)$ is also holomorphic in a neighborhood of $[a_-,a_+]$, the contour integral of this term vanishes. For our purposes, it is not necessary to describe the vector space $\mathrm{Im}\,\mathcal{K}$. Notice that if we apply $\mathcal{K}^{-1}$ to a function $g \in \mathrm{Im}\,\mathcal{K}$ which is furthermore holomorphic outside $\mathbb{C}\setminus [\alpha_-,\alpha_+]$, we can contract the contour $\mathcal{C}([a_-,a_+])$ to a contour $\mathcal{C}([\alpha_-,\alpha_+])$.

\subsubsection{Continuity of $\mathcal{K}$ and $\mathcal{K}^{-1}$}
\label{sec:Kc}
The key fact in this article is that $\mathcal{K}^{-1}$ is a continuous operator in $\big(\mathrm{Im}\,\mathcal{K},\p \cdot \p_{\Gamma}\big)$:
\begin{lemma}
\label{lclo}$\mathrm{Im}\,\mathcal{K}$ is closed subspace of $\mathcal{H}^{(1)}_{1;[a_-,a_+]}$ for the topology induced by the norm $\p \cdot \p_{\Gamma}$, and there exists a constant $k > 0$, such that:
\beq
\forall g \in \mathrm{Im}\,\mathcal{K},\qquad \p \mathcal{K}^{-1}\,g \p_{\Gamma} \leq k \p g \p_{\Gamma}. \nn
\eeq
We call $\p \mathcal{K}^{-1} \p_{\Gamma}$, the infimum of such constants $k$.
\end{lemma}

\noindent \textbf{Proof.} Let us prove first that $\mathcal{K}$, as a endomorphism of $\mathcal{H}_{1;[a_-,a_+]}^{(1)}$, is continuous. For any $f \in \mathcal{H}^{(1)}_{1;[a_-,a_+]}$ in formula~\ref{formuu}, if $x$ runs along $\Gamma$, we can move the contour $\mathcal{C}([a_-,a_+])\cup \mathcal{C}(x)$ to $\Gamma_E$ and get the  bound:
 \bea
 \p \mathcal{K}f \p_{\Gamma} & \leq & 2\big(\|y\|_\Gamma
 \p f \p_{\Gamma} + \frac{\ell(\Gamma_E)}{2\pi}\,
\frac{\p L \p_{\Gamma_E}}{\min_{x\in \Gamma}|L(x)|}\,\Big(\frac{1}{d(\Gamma_E,\Gamma)} + c\Big)\,\p( V^{\{0\}})'\p_{\Gamma_E}\,\p f \p_{\Gamma_E} \nn \\
& \leq & \left[2\,\p y\p_\Gamma + \frac{\ell(\Gamma_E)}{2\pi}\,
\frac{\p L\p_{\Gamma_E}}{\min_{x\in \Gamma}|L(x)|}\Big(\frac{1}{d(\Gamma,\Gamma_E)} + c\Big)\,\p (V^{\{0\}})'\p_{\Gamma_E} \right]\,\p f \p_{\Gamma}. 
\label{conton} 
\eea
$\!\!\!$We have used again the maximum principle for $f$ to find the second line. Likewise, we can show that $\mathcal{K}^{-1}\,:\,\mathrm{Im}\,\mathcal{K} \rightarrow \mathcal{H}_{1;[a_-,a_+]}^{(1)}$ is continuous. In formula~\ref{eq:Km1}, we put $x$ on $\Gamma$, and move the contour from $\mathcal{C}([a_-,a_+])$ to $\Gamma_E$ in Eqn.~\ref{eq:Km1}. Doing so, we pick up a simple pole at $\xi = x$, and we find:
\beq
\big(\mathcal{K}^{-1}g\big)(x) = -\frac{g(x)}{2y(x)} + \frac{1}{\widetilde{\sigma}(x)}\oint_{\Gamma_E} \frac{\dd \xi}{2i\pi}\,\frac{1}{\xi - x}\,\frac{L(\xi)\,g(\xi)}{2S(\xi)}. \nn
\eeq
We deduce :
{\begin{small} \bea
\p \mathcal{K}^{-1} g \p_{\Gamma} & \leq & \frac{\p g \p_{\Gamma}}{2\,\min_{x \in \Gamma} |y(x)|} + \frac{\ell(\Gamma_E)}{4\pi \,d(\Gamma,\Gamma_E)}\,\frac{\max_{\xi \in \Gamma_E} |L(\xi)|}{\min_{x \in \Gamma} |\widetilde{\sigma}(x)|}\,\frac{\p g \p_{\Gamma_E}}{\min_{\xi \in \Gamma_E} |S(\xi)|} \nn \\
& \leq & \left(\frac{1}{2\,\min_{x \in \Gamma} |y(x)|} + \frac{\ell(\Gamma_E)}{4\pi\,d(\Gamma,\Gamma_E)}\,
\frac{\p L\p_{\Gamma_E}}{\min_{\xi \in \Gamma} | \widetilde {\sigma}(x)|
\,\min_{\xi \in \Gamma} |S(\xi)|}\right)\,\p g \p_{\Gamma} . \label{contno} 
\eea
\end{small}}
$\!\!\!$where we used the maximum principle in the last line. Eventually, let us show that $\mathrm{Im}\,\mathcal{K}$ is a closed subspace of $\mathcal{H}_{1;[a_-,a_+]}^{(1)}$. We pick up a sequence $(g_n)_{n}$ in $\mathrm{Im}\,\mathcal{K}$ converging towards $g \in \mathcal{H}^{(1)}_{1;[a_-,a_+]}$ for a norm $\p \cdot \p_{\Gamma_0}$ on a given contour $\Gamma_0$. Let $(f_n)_n$ be a sequence in $\mathcal{H}_{1;[a_-,a_+]}^{(1)}$ such that $g_n = \mathcal{K}f_n$, or equivalently $f_n = \mathcal{K}^{-1}g_n$. Using Eqn.~\ref{contno} for any contour $\Gamma$, we know that $\p f_n \p_{\Gamma} \leq  k\,\p g_n \p_{\Gamma}$ for some constant $k > 0$. So, $f_n$ is a locally bounded subsequence of holomorphic functions in $\mathbb{C}\setminus[a_-,a_+]$. By Montel's theorem, it admits a subsequence $(f_{\phi(n)})_n$ converging to some $f \in \mathcal{H}_{1;[a_-,a_+]}^{(1)}$ uniformly on any compact of $\mathbb{C}\setminus[a_-,a_+]$. Then using Eqn.~\ref{conton}, $g_{\phi(n)} = \mathcal{K}f_{\phi(n)} \rightarrow \mathcal{K}f$ for the norm $\p \cdot \p_{\Gamma_0}$. In particular, $g(x) = \mathcal{K}f(x)$ for all $x \in \mathrm{Ext}(\Gamma_0)$. Since $g$ and $f$ are both analytic in $\mathbb{C}\setminus[a_-,a_+]$, they must coincide on $\mathbb{C}\setminus[a_-,a_+]$. Hence, $g \in \mathrm{Im}\,\mathcal{K}$, showing that $\mathrm{Im}\,\mathcal{K}$ is closed. \hfill $\diamond$

\vspace{0.1cm}

\noindent $\p \mathcal{K}^{-1} \p_{\Gamma}$ is controlled by the distance of the zeroes $s_i$ to the support $[a_-,a_+]$. This motivates Hyp.~\ref{offc}.

\subsubsection{The endomorphism "negative part"}
\label{opN}
Let $g$ be a holomorphic function at least in a neighborhood of $[a_-,a_+]$. The following endomorphism of $\mathcal{H}_{1;[a_-,a_+]}^{(1)}$ often appears in the loop equations:
\beq
\mathcal{N}_{g}(f)(x) = \oint_{\mathcal{C}([a_-,a_+])}\frac{\mathrm{d}\xi}{2i\pi}\,\frac{L(\xi)}{L(x)}\Big(\frac{1}{x - \xi} + c\Big)\,g(\xi)\,f(\xi). \nn
\eeq
We may write sometimes $\mathcal{N}_{g}[f(x)]$ as an abuse of notation. The analyticity assumption on $g$ ensures that $\mathcal{N}_{g}$ is a continuous operator with respect to the norm $\p \cdot \p_{\Gamma}$. Indeed, let us put $x$ on $\Gamma$ and move the contour $\mathcal{C}([a_-,a_+])$ to $\Gamma_E$:
\beq
\mathcal{N}_g(f)(x) = g(x)f(x) + \oint_{\Gamma_E}\frac{\mathrm{d}\xi}{2i\pi}\,\frac{L(\xi)}{L(x)}\Big(\frac{1}{x - \xi} + c\Big)\,g(\xi)\,f(\xi). \nn
\eeq
Thus, the maximum principle implies:
\bea
\p \mathcal{N}_{g}(f) \p_{\Gamma} & \leq &
 \p g \p_{\Gamma} \p f \p_{\Gamma} + \frac{\ell(\Gamma_E)}{2\pi}\Big(\frac{1}{d(\Gamma_E,\Gamma)} + |c|\Big)\,\frac{\mathrm{max}_{\xi \in \Gamma_E} |L(\xi)|}{\mathrm{min}_{x \in \Gamma} |L(x)|}\,\p g \p_{\Gamma_E}\,\p f \p_{\Gamma_E} \nonumber \\
& \leq & \left[\p g \p_{\Gamma} + \frac{\ell(\Gamma_E)} {2\pi}\Big(\frac{1}{d(\Gamma_E,\Gamma)} + |c|\Big)\,\frac{\mathrm{max}_{\xi \in \Gamma_E} |L(\xi)|}{\mathrm{min}_{x \in \Gamma} |L(x)|}\,\p g \p_{\Gamma_E}\right]\p f \p_{\Gamma}. \nn
\eea

\subsection{Order of magnitude of $W_n$}
\label{sec:subs}

If there exists a $1/N$ expansion, $W_n$ ought to be of order of magnitude $N^{2 - n}$. Let us start with a lemma explaining how this can be infered from rough bounds on $W_n$.
Hereafter, $O_l(\cdots)$ or $o_l(\cdots)$ mean $O(\cdots)$ or $o(\cdots)$ with respect to the norm $\p\cdot\p_{\Gamma_l}$. Since the contours $\Gamma_l$ are ordered from the interior to the exterior, being a $o_{l + 1}(\cdots)$ is weaker than being a $o_l(\cdots)$. When the index $l$ is not precised, it is understood that the bound holds for any integer $l$.
\begin{lemma}
\label{lemmaQ}
Let $\delta_{-1}W_1 := N^{-1} W_1 - W_1^{\{-1\}}$ and $l \geq 0$. Assume $\delta_{-1}W_1 \in o_{l}(1)$, and for all integer $n \geq 2$, assume $W_n \in O_{l}(N)$. Then:
\beq
\forall n \geq 2 \qquad\| W_n\|_{\Gamma_{4n - 6 + l}} \in O(N^{2 - n}). \nn
\eeq
\end{lemma}
\textbf{Proof.}  Let $\delta_0 V = V - V^{\{0\}}$. Firstly, as $\delta_{-1}W_1$ and $(\delta_0 V)'$ goes to $0$ uniformly on $\Gamma_{-1}$ when $N \rightarrow \infty$, we observe that  for any fixed integer $k$, and $N$ large enough:
\bea
(1-\varepsilon_{N,k+1})
\p W_n \p_{\Gamma_{k + 1}} & \leq & \p\mathcal{K}^{-1} \p_{\Gamma_{k + 1}}\,\left|\!\left| \Big[ \mathcal{K} + \delta\mathcal{K} + \frac{1}{N}\Big(1 - \frac{2}{\beta}\Big)\frac{\dd}{\dd x}\Big]\,W_n\right|\!\right|_{\Gamma_{k + 1}} \nn \\
\label{eq:gfu}& &  + \frac{1}{N}\Big|1 - \frac{2}{\beta}\Big|\,\zeta_{k}\,\p W_n \p_{\Gamma_{k}},
\eea
where
\bea
[\delta\mathcal{K}](f)(x) & = & -\mathcal{N}_{(\delta_0 V)'}[f(x)] + 2\big(\delta_{-1} W_1\big)(x)\,f(x), \nn \\
\varepsilon_{N,k+1}&=&\p\mathcal{K}^{-1} \p_{\Gamma_{k + 1}}\big(\|{\mathcal N}_{(\delta_0 V)'} \p_{\Gamma_{k + 1}}+ 2\|\delta_{-1} W_1\|_{{\Gamma_{k + 1}}}\big)\eea
goes to zero as $N$ goes to infinity for $k+1\ge l$ by assumption. $\zeta_k$ is defined in Eqn.~\ref{eq:contr}. We assume hereafter that $N$ is large enough so that $\varepsilon_{N,k + 1}$ is smaller than $1/2$.

Secondly, the first version of the loop equation at rank $n \geq 2$ (Thm.~\ref{th:loopeqtn}) can be rewritten:
\beq
\Big[\mathcal{K} + \delta\mathcal{K} + \frac{1}{N}\,\Big(1 - \frac{2}{\beta}\Big)\frac{\dd}{\dd x}\Big]W_{n}(x,x_I) = A_{n + 1} + B_{n} + C_{n - 1} + D_{n - 1},\nn
\eeq
where:
\bea
A_{n + 1} & = & -\frac{1}{N}\,W_{n + 1}(x,x,x_I), \nn \\
B_{n} & = & - \frac{1}{N}\sum_{\substack{n_1,n_2 \geq 1 \\ n_1 + n_2 = n - 1}} \sum_{\substack{J \subseteq I \\ |J| = n_1}} W_{n_1 + 1}(x,x_J)\,W_{n_2 + 1}(x,x_{I\setminus J}), \nn \\
C_{n - 1} & = & -\frac{1}{N}\,\frac{2}{\beta}\sum_{i \in I} \frac{\dd}{\dd x_i}\left\{\frac{W_{n - 1}(x,x_{I\setminus\{i\}})}{x - x_i} - \frac{L(x_i)}{L(x)}\Big(\frac{1}{x - x_i} + c\Big) W_{n - 1}(x_I)\right\}, \nn \\
D_{n - 1} & = & \frac{1}{N}\,\frac{2}{\beta}\sum_{\tau\in \mathsf{Soft}} \frac{\partial_{a_\tau} W_{n - 1}(x_I)}{x - a_\tau}. \nn
\eea
We know from Proposition~\ref{pro} that $D_{n - 1} \in O(e^{-N\,\eta})$, so this term does not contribute at any order of magnitude $N^{-k}$.
Now, if we assume that $W_n \in O_{l}(N)$ for all $n \geq 2$ (this is obviously true for $n = 1$), we always have $A_{n + 1} \in O_{l}(1)$ and $C_{n - 1} \in O_{l + 2}(1)$, whereas the last  last term in Eqn.~\ref{eq:gfu} is bounded by hypothesis for $k\ge l$.

\begin{figure}
\begin{center}
\includegraphics[width=0.8\textwidth]{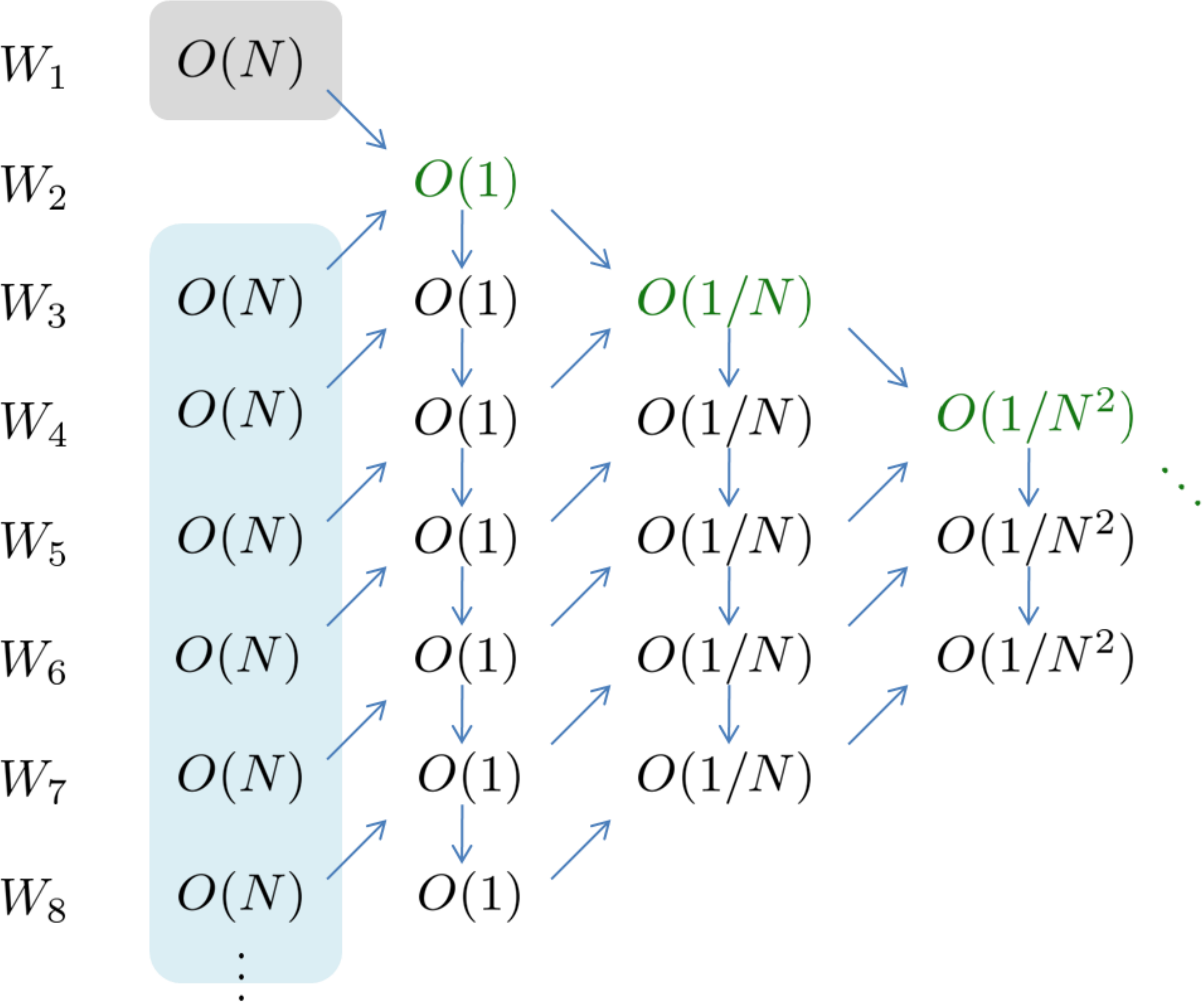}
\caption{\label{fig:reco} Triangular recursion for Lemma~\ref{lemmaQ}}
\end{center}
\end{figure}

Now, we want to bound $W_n$ by induction on $n$. At rank $n =2$, we have $B_2 = 0$, and we deduce from Eqn.~\ref{eq:gfu} that $W_2 \in O_{l + 2}(1)$. Then at rank $n = 3$, the product term $B_3$ is $O_{l + 2}(1/N)$ and $C_2$ is $O_{l + 4}(1/N)$,  thus $W_3 \in O_{l + 4}(A_4) = O_6(1)$. Then similarly at rank $n = 4$, the product term $B_4$ is $O_{l + 4}(1/N)$ and $C_3$ is $O_{l + 6}(1/N)$, thus $W_4 \in O_{l + 6}(1)$. This implies in return that $A_4 \in O_{l + 6}(1/N)$, thus $W_3 \in O_{l + 6}(1/N)$. And so on \ldots{} The result can be proved by a triangular induction, as depicted in Fig.~\ref{fig:reco}. At each vertical step, we are forced to trade the contour $\Gamma_k$ with the exterior contour $\Gamma_{k + 2}$ in order to control the $C$ terms. So, to go from $W_n \in O_{k_n}(N^{2 - n})$ (in the $n^{\mathrm{th}}$ column) to $W_n \in O_{k_{n + 1}}(N^{2 - (n + 1)})$ (in the $(n + 1)^{\mathrm{th}}$ column), we must reach $W_{n + 2}$ in the $\mathrm{n}^{\mathrm{th}}$ column. This is done by two vertical steps, thus $k_{n + 1} = k_{n} + 4$. Since $k_2 = l + 2$, we have $k_n = 4n - 6 + l$ for all $n \geq 2$.  \hfill $\diamond$

\begin{lemma}
\label{lemmaR}
If there exists $\gamma \in [0,1[$ and $\delta \in [0,\infty[$ such that $W_n \in O_0(N^{\gamma n - \delta})$ for all $n \geq 2$, then:
\beq
\p W_n \p_{\Gamma_{4n - 6 + l}} \in O(N^{2 - n}), \nn
\eeq
where $l = 2 \lceil (\gamma^{-1} - 1)^{-1} \rceil$.
\end{lemma}
\textbf{Proof.} Now, let us rather assume the existence of $\gamma \in [0,1[$ and $\delta \geq 0$ such that, for all $n \geq 2$, $W_n \in O_0(N^{\gamma n - \delta})$. $D_{n - 1}$ being always exponentially small, it does not matter in our discussion. At rank $n = 2$, as $B_{2} = 0$ and $C_{1} \in O_2(1)$, we have $W_2 \in O_2\big(\mathrm{max}[\frac{1}{N}\,W_3,1]\big)$. We also have for all $n$:
\bea
A_{n + 1} \in O_0\big(\mathrm{max}[N^{\gamma n - \delta - (1 -  \gamma)},1]\big), && \nn \\
B_n \in O_0\big(\mathrm{max}[N^{\gamma n - 2\delta - (1 - \gamma)},1]\big), && \nn \\
C_{n - 1} \in O_2\big(\mathrm{max}[N^{\gamma n - \delta - (1 + \gamma)},1]\big). \nn
\eea
When these $O(\cdots)$ decay, it does not hurt to consider them as $O(1)$. So, our bounds are upgraded at least to $W_n \in O_2\big(\mathrm{max}[N^{\gamma n - \delta'},1]\big)$ with $\delta' = \delta + 1 - \gamma > \delta$. By repeating the argument $k$ times, we obtain for all $n \geq 2$, $W_n \in O_{2k}\big(\mathrm{max}[N^{((k + 1)\gamma - k)n - \delta},1]\big)$. Since $\gamma < 1$, by choosing an integer $k \geq \frac{1}{1/\gamma - 1}$, we deduce that $W_n \in O_{l}(1)$ for all $n \geq 2$ with $l = 2k$, and we apply Lemma~\ref{lemmaQ} to conclude. \hfill $\diamond$

\subsection{Initialization}
\label{secini}
We now establish a priori control on the correlators. We shall use:

\begin{lemma}\label{sec:assum}
Let $w_N=N^\epsilon$ for some $\epsilon>0$.
Assume that
 for any integer $p$,
there exists $C_p > 0$ and independent of $N$, such that for all $x\in \mathbb{C}\setminus[a_-,a_+]$:
\beq
\label{eq:27}\mu_{N,\beta}^V\Big\{\Big|\int \frac{\dd M_N(\xi)}{x - \xi} -
\mu_{N,\beta}^V\Big[\int \frac{\dd M_N(\xi)}{x - \xi}\Big]\Big|^p\Big\} \le \frac{C_p w_N^p}{\big(d(x,[a_-,a_+])\big)^{2p}}.
\eeq
Then, for all $n \geq 2$, $W_n \in O(w_N^{n})$ for the norm $\p\cdot\p_{\Gamma}$, when $N \rightarrow \infty$.
\end{lemma}
\noindent {\bf Proof}. For $n \ge 2$, $W_n(x_1,\ldots,x_n)$ is a polynomial in:
\beq
\mu_{N,\beta}^{V;[a_-,a_+]}\Big\{\prod_{j \in J}
\Big(\int \frac{\dd M_N(\xi)}{x_j - \xi} - \mu_{N,\beta}^V\Big[\int \frac{\dd M_N(\xi)}{x_j - \xi}\Big]\Big)\Big\}, \nn
\eeq
with $J \subseteq \{1,\ldots,n\}$, and the coefficients of this polynomial are independent of
$N$. Thus by Eqn.~\ref{eq:27} and H\"older inequality, there exists $D_n \in \mathbb{R}_+^*$ independent of $N$ such that:
\beq
\big|W_n(x_1,\ldots,x_n)\big|\le \frac{D_n\,w_N^n}{\big(\min_{1\le i\le n} d(x_i,[a_-,a_+])\big)^{2n}}. \nn
\eeq
Hence, taking the sup for $x_i \in \Gamma$, we find $W_n \in O(w_N^n)$.

\begin{lemma}
\label{sec:assumconv}
Under the five assumptions of Section~\ref{hypotheses}, Eqn.~\ref{eq:27} holds for any $\epsilon > 0$.
\end{lemma}
\textbf{Proof.} Our starting point comes from a result of Boutet de
 Monvel, Pastur and Shcherbina \cite{Boutet}, developed
 by Johansson\footnote{Johansson's has written his proof in the framework $[a_-,a_+] = \mathbb{R}$, but there is no difficulty adapting it to $[a_-,a_+]$ finite. $\mathrm{Im}\,z$ should be replaced by $d(z,[a_-,a_+])$, and its powers in the bound of his Lemma 3.10 and 3.11 may differ, but the order of magnitude $\omega_N$ (our $w_N$) is the same.}\cite[(3.49)]{Johansson} and more recently in \cite[(2.26)]{Shch2}. Let $\Gamma' \subseteq \mathrm{Int}\,\Gamma$ be a contour surronding $[a_-,a_+]$. For any $\phi\,:\,\mathrm{Int}(\Gamma) \rightarrow \mathbb{C}$ which is a continuous function, and real-valued on $[a_-,a_+]$, there exists a positive constant $C$ such that:
\beq
\mu_{N,\beta}^{V}\left[\exp\left(\frac{1}{2\big(\sup_{z' \in \Gamma'} |\phi(z')|\big)w_N}\Big(\int \phi(\xi)\dd M_N(\xi) - N\int \phi(\xi)\dd L(\xi)\Big)
\right) \right] \leq 3, \nn
\eeq
where $w_N = C\ln N$. By Chebychev's inequality, we deduce that:
\beq
\forall t \in [0,+\infty[,\quad \mu_{N,\beta}^{V}\left\{\Big|\int \phi(\xi)\dd M_N(\xi) - N\int\phi(\xi)\dd L(\xi)\Big| \geq t\big(\sup_{z' \in \Gamma'} |\phi(z')|\big)w_N\right\} \leq 6e^{-t}, \nn
\eeq
and therefore for all $p \in \mathbb{N}$, 
\beq
 \mu_{N,\beta}^{V}\left[\Big|\int \phi(\xi)\dd M_N(\xi) - N\int \phi(\xi)\dd L(\xi)\Big|^{p}\right] \leq p!\big(\mathrm{sup}_{z' \in \Gamma'} |\phi(z')|\big)^{p}w_N^p. \nn
\eeq
In particular, we can apply this discussion to $\phi(z) = \mathrm{Re}\,\frac{1}{x - z}$ and $\phi(z) = \mathrm{Im}\,\frac{1}{x - z}$ where $x$ is a point of $\Gamma$. This leads to Eqn.~\ref{eq:27}. \hfill $\diamond$

In the case of a strictly convex potential, we may  use instead concentration of measure:
\begin{lemma}
If $V^{\{0\}}$ is strictly convex on $[a_-,a_+]$, then Eqn.~\ref{eq:27} holds with $\e=0$.
\end{lemma}
{\bf Proof.}
Since $V^{\{0\}}$ is strictly convex on $[a_-,a_+]$, $V$ is also strictly convex on $[a_-,a_+]$ for $N$ large enough. By concentration of measure, see \cite{GZ} or \cite[Section 2.3 and Exercise 4.4.33]{AGZ}, there exists $c > 0$
such that, for all $x\in\mathbb C\backslash \mathbb [a_-,a_+]$, for all $\epsilon>0$
and $N\in \mathbb N$:
\beq
\mu_{N,\beta}^V\left\{\Big|\int \frac{\dd M_N(\xi)}{x - \xi} - \mu_{N,\beta}^V\Big[\int \frac{\dd M_N(\xi)}{x - \xi}\Big]\Big|\ge \frac{\epsilon}{\big(d(x,[a_-,a_+])\big)^2} \right\} \le 2e^{-c\epsilon^2}. \nn
\eeq
This entails Eqn.~\ref{eq:27}. \hfill $\diamond$

\subsection{Leading order of $W_1$}
\label{sec:kos}
Afterwards, all steps only rely on the analysis of loop equations. Although we already know the characterization of the equilibrium measure $\mu_{\mathrm{eq}}$, and thus of its Stieltjes transform $W_1^{\{-1\}}$, let us recall how $W_1^{\{-1\}}$ is characterized by the loop equations. We write the loop equation at rank $1$ (Thm.~\ref{th:loopeqt}):

{\small \bea
 \frac{1}{N^2}\,W_2(x,x) &&  \nn \\
+ \big(W_1^{\{-1\}}(x)\big)^2  - \frac{1}{(x - a_-)(x - a_+)} && \nn \\
 - \oint_{{\mathcal C}([a_-,a_+])} \frac{\dd\xi}{2i\pi}\,\frac{1}{x - \xi}\,\frac{(\xi - a_-)(\xi - a_+)}{(x - a_-)(x - a_+)}\,(V^{\{0\}})'(\xi)\,W_1^{\{-1\}}(\xi) && \nn \\
+ \mathcal{K}\big[\delta_{-1} W_1\big](x) + \frac{1}{N}\,\Big(1 - \frac{2}{\beta}\Big)\Big(W_1^{\{-1\}}(x) + \frac{1}{(x - a_-)(x - a_+)}\Big) - \frac{1}{N}\,\mathcal{N}_{V^{\{1\}}}\big[W_1^{\{-1\}}\big](x) && \nn \\
 + \big((\delta_{-1}W_1)(x)\big)^2  - \mathcal{N}_{(\delta_0 V)'}\big[\delta_{-1}W_1](x) + \frac{1}{N}\,\mathcal{N}_{(\delta_1 V)'}\big[W_1^{\{-1\}}\big](x) & = & 0. \nn \\
\label{theeq2} &&
\eea}
$\!\!\!$We already know that the $4^{\mathrm{th}}$ and the $5^{\mathrm{th}}$ line are $o(1)$. Since $W_2 \in o(N^2)$, $W_1^{\{-1\}}$ satisfy the loop equation at leading order:

{\small \bea
\label{theeq}  \big(W_1^{\{-1\}}(x)\big)^2 & =& \frac{1}{(x - a_-)(x - a_+)} \\
&& + \oint_{{\mathcal C}([a_-,a_+])} \frac{\dd\xi}{2i\pi}\,\frac{1}{x - \xi}\,\frac{(\xi - a_-)(\xi - a_+)}{(x - a_-)(x - a_+)}\,(V^{\{0\}})'(\xi)\,W_1^{\{-1\}}(\xi). \nn
\eea}

\begin{remark}
\label{Rem}
Recall that $\mathrm{supp}\,\rho = [\alpha_-,\alpha_+]$ is the discontinuity locus of $W_1^{\{-1\}}$. By the properties of the Stieltjes transform:
\beq
y(x) = \frac{(V^{\{0\}})'(x)}{2} - W_1^{\{-1\}}(x) \nn
\eeq
defines a holomorphic function on $U\setminus[\alpha_-,\alpha_+]$, and:
\beq
\forall x_0 \in [\alpha_-,\alpha_+],\qquad \lim_{\epsilon \rightarrow 0^+} y(x_0 + i\epsilon) = i\pi\,\rho(x_0).\nn
\eeq
We state that there exists $M(x)$, continuous in some open neighborhood of $[\alpha_-,\alpha_+]$, such that:
\beq
\label{eq:remsk1} y(x) = \frac{M(x)}{\sqrt{(x - \alpha_-)(x - \alpha_+)}}.
\eeq
\end{remark}
\textbf{Proof.} In Eqn.~\ref{theeq}, we may first deform the contour $\mathcal{C}([a_-,a_+])$ to $\mathcal{C}([\alpha_-,\alpha_+])$. Secondly, we can rewrite:
\bea
&& \big(W_1^{\{-1\}}(x)\big)^2 - (V^{\{0\}})'(x)W_1^{\{-1\}}(x)  + \frac{U(x)}{(x - a_-)(x - a_+)} = 0, \nn \\
&& U(x) = -1 + \oint_{\mathcal{C}([\alpha_-,\alpha_+])\cup\mathcal{C}(x)}\frac{\dd\xi}{2i\pi}\,\frac{(\xi - a_-)(\xi - a_+)}{x - \xi}\,(V^{\{0\}})'(\xi)\,W_1^{\{-1\}}(\xi),  \nn
\eea
where $U(x)$ is now holomorphic in some open neighborhood of $[\alpha_-,\alpha_+]$. So:
\beq
y(x) = \sqrt{\frac{R(x)}{(x - a_-)(x - a_+)}},\quad R(x) = \frac{1}{4}(x - a_-)(x - a_+)(V^{\{0\}})'(x) - U(x). \nn
\eeq
This equation tells us that  the discontinuity of $y$ is of squareroot type. If $\alpha_- = a_-$ and $\alpha_+ = a_+$, we have Eqn.~\ref{eq:remsk1}.
If say $a_- < \alpha_-$, the fact that $y(x)$ has no discontinuity on $[a_-,\alpha_-[$ but
a discontinuity on $[\alpha_-,\alpha_+]$ forces $R(x)$ to have a simple zero at $x = a_-$ and at $x = \alpha_-$, so that $y(x)$ is finite when $x = a_-$ and vanishes as $O(\sqrt{x - \alpha_-})$ when $x \rightarrow \alpha_-$. A similar statement holds if $a_+ > \alpha_+$. Then, Eqn.~\ref{eq:remsk1} holds a fortiori. \hfill $\diamond$

\subsection{First correction to $W_1$}

Let us reconsider Eqn.~\ref{theeq2} (or the equivalent relation taking Remark~\ref{sechy} into account) after removing the $2^{\mathrm{nd}}$ and the $3^{\mathrm{rd}}$ line which has just been identified as the leading order. We can write as in \S~\ref{sec:subs}:
\beq
\left[\mathcal{K} + \widetilde{\delta\mathcal{K}} + \frac{1}{N}\Big(1 - \frac{2}{\beta}\Big)\frac{\dd}{\dd x}\right]\delta_{-1}W_1(x) = A_2 + C_0 + D_0, \nn
\eeq
where:
\bea
\widetilde{\delta\mathcal{K}}[f](x) & = & -\mathcal{N}_{(\delta_0V)'}[f(x)] + \delta_{-1}W_1(x)f(x), \nn \\
A_2 & = & -\frac{1}{N^2} W_2(x,x), \nn \\
C_0 & = & -\frac{1}{N}\Big(1 - \frac{2}{\beta}\Big)\left(\sum_{\tau\in\mathsf{Hard}} \frac{1}{a_\tau - a_{-\tau}}\,\frac{1}{x - a_\tau}\right), \nn \\
D_0 & = & \sum_{\tau\in\mathsf{Soft}} \frac{\partial_{a_\tau} \ln Z}{x - a_\tau}. \nn
\eea
By an argument similar to Eqn.~\ref{eq:gfu}, knowing that $W_2 \in O_l(N)$ implies that $\delta_{-1}W_1 \in O_{l + 1}(1/N)$. Assuming further $W_3 \in O_{l'}(N)$ implies after Section~\ref{sec:subs} that $W_2 \in O_{l' + 2}(1)$, so the $1^{\mathrm{st}}$ line of Eqn.~\ref{theeq2} is subleading compared to the $3^{\mathrm{rd}}$ line. These two bounds are provided by Section~\ref{secini} (the values of $l$ and $l'$ do not matter here). Hence:

\begin{lemma}
\label{lema} There exists $W_1^{\{0\}} \in \mathcal{H}_{1;[\alpha_-,\alpha_+]}^{(1)}$ such that $W_1 = NW_1^{\{-1\}} + W_1^{\{0\}} + o(1)$. Explicitly:
{\small\beq
W_1^{\{0\}}(x) = \mathcal{K}^{-1}\left\{-\Big(1 - \frac{2}{\beta}\Big)\left[\frac{\dd}{\dd x}\big(W_1^{\{-1\}}(x)\big) + \sum_{\tau \in \mathsf{Hard}} \frac{1}{a_\tau - a_{-\tau}}\,\frac{1}{x - a_{\tau}}\right] +  \mathcal{N}_{(V^{\{1\}})'}(W_1^{\{-1\}})(x)\right\}.\nn
\eeq}
\end{lemma}
This order was also obtained by \cite{Shch2} with similar arguments.

\subsection{Recursion hypothesis at order $k_0$}
\label{sec:recc}
Let $k_0 \geq -1$. We assume that the correlators $W_n$ (for all $n \geq 1$) are determined up to a $o(N^{-k_0})$ for the norm $\p \cdot \p_{\Gamma_{l(k_0;n)}}$.
\beq
\label{eq:expk0}W_{n}(x_1,\ldots,x_{n}) = \sum_{k = n - 2}^{k_0} N^{-k}\,W_{n}^{\{k\}}(x_1,\ldots,x_n) + N^{-k_0}\,\delta_{k_0}W_n(x_1,\ldots,x_n).
\eeq
Here, $W_n^{\{k\}}(x_1,\ldots,x_n)$ are already known (they depend on $\beta$ but not on $N$), and we call:
\beq
\omega_{n}^{\{k\}} = \sup_{\substack{1 \leq n' \leq n \\ -1 \leq k' \leq k}} \p W_n^{\{k\}} \p_{\Gamma_{l(k_0;n)}}, \nn
\eeq
a bound for their norm. We can always assume that $l(k,n)$ defined for $-1 \leq k \leq k_0$ and $n \geq 1$ is an increasing function of $k$ and $n$. Though the errors $\delta_{k_0}W_n$ are not supposed to be known, we assume that they are small:
\beq
\forall n \geq 1 \qquad \p \delta_{k_0} W_{n} \p_{\Gamma_{l(k_0;n)}} \leq \epsilon_{N}^{\{k_0\}}\,\Delta_{n}^{\{k_0\}}. \nn
\eeq
Here, $\epsilon_{N}^{\{k_0\}}$ depends only on $N$ and $k_0$, and $\epsilon_{N}^{\{k_0\}} \rightarrow 0$ when $N \rightarrow \infty$, and $\Delta_{n}^{\{k_0\}}$ is a constant independent of $N$. We may assume that $\Delta_{n}^{\{k_0\}}$ increases with $n \geq 1$, upon replacement by $\sup_{1 \leq n' \leq n} \Delta_{n}^{\{k_0\}}$. When $n > k_0 + 2$, we assume that Eqn.~\ref{eq:expk0} reduces to:
\beq
W_{n} = N^{-k_0}\,\delta_{k_0}W_{n}. \nn
\eeq
Lemma~\ref{sec:assum} and Section~\ref{lema} ensure that the initial ($k_0 = -1$) recursion hypothesis is satisfied. Moreover, we can take $\epsilon_{N}^{\{-1\}} = 1/N$, and up to a redefinition $\Gamma_k \rightarrow \Gamma_{k - m}$ for some integer $m$, we can take $l(-1;n) = 4(n - 1)$.

\subsection{Determination of $\delta_{k_0}W_{n_0}$}
\label{sec:bboun}
Let $n_0 \geq 1$. We now turn to the determination of the leading order of $\delta_{k_0}W_{n_0}(x,x_I)$. The case $(n_0,k_0) = (1,-1)$ is a bit special (because of the second term of the second line in Eqn.~\ref{theeq2}) and is given by Lemma~\ref{lema}. In all other cases, we consider the loop equation at rank $n_0$ (Thm.~\ref{th:loopeqtn}). Up to $o(N^{-(k_0 - 1)})$, the equation is true and involves quantities which are already known from the recursion hypothesis. The equality of the $o(N^{-(k_0 - 1)})$ involves the unknown $\delta_{k_0}W_{n_0}(x,x_I)$. The operator $\mathcal{K}$ introduced in \S~\ref{opK} plays a special role. When the potential $V$ has a $1/N$ expansion, the operator $\mathcal{N}$ introduced in \S~\ref{opN} also appears, and we denote:
\beq
V = \sum_{k = 0}^{k_0 + 1} N^{-k}\,V^{\{k\}} + N^{-(k_0 + 1)}\,\delta_{k_0 + 1} V. \nonumber
\eeq
We find:
\beq
\label{equareste} N^{-(k_0 - 1)}\,\mathcal{K}\big(\delta_{k_0}W_{n_0}\big)(x,x_I) = - N^{-k_0}\,E^{\{k_0\}}_{n_0}(x,x_I) - N^{-k_0}\,R^{\{k_0\}}_{n_0}(x,x_I),
\eeq
with
\bea
E_{n_0}^{\{k_0\}}(x,x_I) & := & W_{n_0 + 1}^{\{k_0\}}(x,x,x_I)  - \sum_{k = 1}^{k_0 + 1} \mathcal{N}_{(V^{\{k\}})'}\big[W_{n_0}^{\{k_0 + 1 - k\}}(x,x_I)\big]  \nn \\
& & + \sum_{J \subseteq I} \sum_{k = 0}^{k_0} W_{|J| + 1}^{\{k\}}(x,x_J)\,W_{n_0 - |J|}^{\{k_0 - k\}}(x,x_{I\setminus J})+ \Big(1 - \frac{2}{\beta}\Big)\frac{\mathrm{d}}{\mathrm{d} x}\Big(W_{n_0}^{\{k_0\}}(x,x_I)\Big) \nonumber \\
&& + \frac{2}{\beta}\,\sum_{i \in I} \frac{\mathrm{d}}{\mathrm{d} x_i}\left\{\frac{W_{n_0 - 1}^{\{k_0\}}(x,x_{I\setminus\{i\}})}{x - x_i} - \frac{L(x_i)}{L(x)}\Big(\frac{1}{x - x_i} + c\Big)W_{n_0 - 1}^{\{k_0\}}(x_I)\right\}, \nonumber
\eea
and the remaining
\bea
R^{\{k_0\}}_{n_0}(x,x_I)& := &  \delta_{k_0}W_{n_0 + 1}(x,x,x_I)  + \sum_{k = 1}^{k_0} N^{-k}  \sum_{k' = 0}^{k_0} \sum_{J \subseteq I} W_{|J| + 1}^{\{k'\}}(x,x_J)\,W_{n_0 - |J|}^{\{k_0 + k - k'\}}(x,x_{I \setminus J}) \nonumber \\
&& + \sum_{k = 0}^{k_0} N^{-k} \sum_{J \subseteq I} \big(\delta_{k_0}W_{|J| + 1}\big)(x,x_{J})\,W_{n_0 - |J|}^{\{k\}}(x,x_{I\setminus J}) \nonumber \\
& & + N^{-k_0}\sum_{J \subseteq I} \big(\delta_{k_0}W_{|J| + 1}\big)(x,x_J)\,\big(\delta_{k_0}W_{n_0 - |J|}\big)(x,x_{I \setminus J}) \nonumber \\
& & + \Big(1 - \frac{2}{\beta}\Big)\frac{\mathrm{d}}{\mathrm{d}x}\Big(\big(\delta_{k_0}W_{n_0}\big)(x,x_I)\Big)- \big(\delta_{k_0}W_{n_0 + 1}\big)(x,x,x_I) \nonumber \\
&&  - \sum_{k = 0}^{k_0} N^{-k}\,\mathcal{N}_{\big(\delta_{k_0 + 1}V\big)'}\:\big[W_{n_0}^{\{k\}}(x,x_I)\big]  - \sum_{k = 0}^{k_0} N^{-k}\,\mathcal{N}_{(V^{\{k + 1\}})'}\big[\big(\delta_{k_0}W_{n_0}\big)(x,x_I)\big] \nonumber \\
&& - N^{-k_0}\,\mathcal{N}_{\big(\delta_{k_0 + 1}V\big)'}\big[\big(\delta_{k_0}W_{n_0}\big)(x,x_I)\big] \nn \\
&& + \frac{2}{\beta}\sum_{i \in I}\frac{\mathrm{d}}{\mathrm{d}x_i}\left\{\frac{\big(\delta_{k_0}W_{n_0 - 1}\big)(x,x_{I\setminus\{i\}})}{x - x_i} - \frac{L(x_i)}{L(x)}\Big(\frac{1}{x - x_i} + c\Big)\big(\delta_{k_0}W_{n_0 - 1}\big)(x_I)\right\} \nonumber \\
\label{eq:error} && + \frac{2}{\beta}  \sum_{\tau \in \mathsf{Soft}} N^{k_0}\,\frac{\partial_{a_\tau} W_{n - 1}(x_I)}{x - a_\tau}. 
\eea
It is understood that $\mathcal{K}$ and $\mathcal{N}_g$ operate on the $x$ variable. The variables $x_I$ are spectators.
Notice that this equation is linear in $\delta_{k_0}W_{n_0}$, up to a small quadratic term.

Looking naively at this equation, we see that the leading term of $\delta_{k_0}W_{n_0}$ happens to be of order $1/N$ (giving a $N^{-(k_0 + 1)}$ contribution to $W_{n_0}$), and is obtained by applying $\mathcal{K}^{-1}$ to $E_{n_0}^{\{k_0\}}(x,x_I)$. To make this idea rigorous, let us bound $R_{n_0}^{\{k_0\}}$. Even if some terms in the right hand side have not been determined yet (like $\delta_{k_0}W_{n}$ that we are just considering), we already know a bound for each of them from the recursion hypothesis. Very rough bounds are enough, we just need to show that the right hand side is small when $N \rightarrow \infty$. When $k_0 = -1$, we must pay special attention at the terms involving $N^{-k_0}$ directly, i.e. the $3^{\textrm{rd}}$ line and the $6^{\textrm{th}}$ line in Eqn.~\ref{eq:error}. In the $6^{\textrm{th}}$ line, $(\delta_{k_0 + 1}V)'$ is of order $N^{-1}$, so we obtain a term of order $\epsilon_{N}^{\{k_0\}}$, which is always small. The $3^{\mathrm{rd}}$ line is of order $N(\epsilon_{N}^{\{-1\}})^2$, which is also small since we have here $\epsilon_{N}^{\{-1\}} = 1/N$ (Lemma~\ref{lema}). For $N$ large enough, we have:
\bea
\p R_{n_0}^{\{k_0\}} \p_{\Gamma_{l(k_0 + 1;n_0)}} & \leq & \epsilon_{N}^{\{k_0\}}\,\Delta_{n_0 + 1}^{\{k_0\}}  + N^{-1}\,(k_0 + 1)2^{n_0 - 1}\,\big(\omega_{n_0}^{\{k_0\}}\big)^2 \nonumber \\
& & + \epsilon_{N}^{\{k_0\}}\,2^{n_0 - 1}\,\Delta_{n_0}^{\{k_0\}}\,\omega_{n_0}^{\{k_0\}} + (\epsilon_{N}^{\{k_0\}})^2\,N^{-k_0}\,\big(\Delta_{n_0}^{\{k_0\}}\big)^2 \nonumber \\
& & + \epsilon_{N}^{\{k_0\}}\,\Big|1 - \frac{2}{\beta}\Big|\,\zeta_{l(k_0;n_0)}\,\Delta_{n_0}^{\{k_0\}} + N^{-1}\,\sum_{k = 0}^{k_0} \p\mathcal{N}_{\big(\delta_{k_0 + 1}V)'}\p_{\Gamma_{l(k_0;n_0)}}\,\Delta_{n_0}^{\{k_0\}} \nonumber \\
& & + \epsilon_{N}^{\{k_0\}}\,\sum_{k = 0}^{k_0} \p \mathcal{N}_{(V^{\{k\}})'}\p_{\Gamma_{l(k_0;n_0)}}\,\Delta_{n_0}^{\{k_0\}} + \epsilon_{N}^{\{k_0\}}\,N^{-k_0}\,\p \mathcal{N}_{\big(\delta_{k_0 + 1}V\big)'}\p_{\Gamma_{l(k_0;n_0)}}\,\Delta_{n_0}^{\{k_0\}} \nonumber \\
& & + \epsilon_{N}^{\{k_0\}}\,\frac{2}{\beta}\,\zeta_{l(k_0;n_0 - 1)}\,\left(|c| + \zeta_{l(k_0;n_0 - 1) + 1}\,\frac{\mathrm{sup}_{\xi \in \Gamma_{l(k_0;n_0 - 1)}} |L(\xi)|}{\inf_{x \in \Gamma_{l(k_0;n_0 - 1)}} |L(x)|}\right)\,\Delta_{n_0 - 1}^{\{k_0\}} \nn \\
&& + \frac{2}{\beta}\,\frac{2\gamma_n}{d(\Gamma,[a_-,b_-])^{n}} N^{k_0 + n - 1}\,e^{-N\,\eta} \nn
\eea
Given the control provided by the recursion hypothesis, this inequality is correct provided we choose:
\beq
\label{eq:order}l(k_0 + 1;n_0) \geq \mathrm{max}\big[l(k_0;n_0 - 1) + 2,l(k_0;n_0) + 1,l(k_0;n_0 + 1)\big].
\eeq
Accordingly, $R_{n_0}^{\{k_0\}} \rightarrow 0$ when $N \rightarrow \infty$. Eqn.~\ref{equareste} tells us that $E_{n_0}^{\{k_0\}} + R_{n_0}^{\{k_0\}} \in \mathrm{Im}\,\mathcal{K}$ for any $N$. Since $\mathrm{Im}\,\mathcal{K}$ is closed (Lemma~\ref{lclo}), we know that $E_{n_0}^{\{k_0\}} \in \mathrm{Im}\,\mathcal{K}$, and also by difference $R_{n_0}^{\{k_0\}} \in \mathrm{Im}\,\mathcal{K}$ for any $N$. And, by continuity of $\mathcal{K}^{-1}$, we deduce:
\beq
\delta_{k_0}W_{n_0} = \frac{1}{N}\,W_{n_0}^{\{k_0 + 1\}} + \frac{1}{N}\,\delta_{k_0 + 1}W_{n_0}, \nn
\eeq
where:
\beq
\label{comi}W_{n_0}^{\{k_0 + 1\}} = -\mathcal{K}^{-1}[E_{n_0}^{\{k_0\}}],\qquad \delta_{k_0 + 1}W_{n_0} = -\mathcal{K}^{-1}[R_{n_0}^{\{k_0\}}] \in o(1).
\eeq
The previous inequality is more precise about the $o_{l(k_0 + 1;n_0)}(1)$: there exists a constant $\Delta_{n_0}^{\{k_0 + 1\}}$, such that
\beq
\p \delta_{k_0 + 1} W_{n_0} \p_{\Gamma_{l(k_0 + 1;n_0)}} \leq \Delta_{n_0}^{\{k_0 + 1\}}\,\mathrm{max}(N^{-1}\,;\,\epsilon_{N}^{\{k_0\}}). \nn
\eeq

\subsection{Remarks}

The recursion hypothesis tells us that $W_{n_0}^{\{k_0\}} = 0$ whenever $n > k_0 + 2$ (we call $\star[k_0]$ this recursive assumption). Let us see what happens at order $k_0 + 1$ (here, $k_0$ is fixed, but $n_0$ is free), by looking at Eqn.~\ref{comi}.
\begin{itemize}
\item[$\bullet$] The term $W_{n_0 + 1}^{\{k_0\}}$ vanishes whenever $n_0 > k_0 + 1$.
\item[$\bullet$] The term $W_{|J| + 1}^{\{k\}}\,W_{n_0 - |J|}^{\{k_0 - k\}}$ may be non zero in case $k + 1 \geq |J| \geq n_0 - k_0 - 2 + k$. This is impossible to fulfil as soon as $n_0 > k_0 + 3$.
\item[$\bullet$] The term $\Big(W_{n_0}^{\{k_0\}}\Big)'$ vanishes whenever $n_0 > k_0 + 2$.
\item[$\bullet$] The term involving $W_{n_0 - 1}^{\{k_0\}}$ vanishes whenever $n_0 > k_0 + 3$.
\end{itemize}
Accordingly, $W_{n_0}^{\{k_0 + 1\}} \equiv 0$ when $n_0 > (k_0 + 1) + 2$, i.e. $\star[k_0 + 1]$ holds. This is just the manifestation of Lemma~\ref{lemmaQ}.
Hence, we have propagated the full recursion hypothesis to order $k_0 + 1$. An easy recursion shows that $W_n^{\{k\}}$ are actually holomorphic functions on the domain $\mathbb{C}\setminus[\alpha_-,\alpha_+]$, i.e belongs to the subspace $\mathcal{H}_{n;[\alpha_-,\alpha_+]}^{(1)}$ of $\mathcal{H}_{n;[a_-,a_+]}^{(1)}$. Therefore, we can contract the contour to $\mathcal{C}([\alpha_-,\alpha_+])$ in the expression of $\mathcal{K}^{-1}$ (Eqn.~\ref{eq:Km1}) when computing $W_n^{\{k\}}$ with formula~\ref{comi}.

Since $l(k;n) = 4(n - 1)$, the minimal solution of Eqn.~\ref{eq:order} is $l(k;n) = 4(n + k)$. Indeed, in this proof, we need to have a more restrictive control on the error done at height $n + k$, in order to bound the error done at height $n + k + 1$. Nevertheless, since $\Gamma_l \subseteq \mathrm{Int}(\Gamma_E)$ for all $l$, we can at the end make the weaker statement that, for any $n$ and $k$:
\beq
\label{errcontrol}\p \delta_{k}W_n \p_{\Gamma_E} \rightarrow 0
\eeq
when $N \rightarrow \infty$. However, we necessarily have
$d(\Gamma_l,\Gamma_{l + 1}) \rightarrow 0$ when $l \rightarrow \infty$, so that the constant $\zeta_{l}$ which allows us to bound the derivative of a function with the function itself (Eqn.~\ref{eq:contr}), blows up. This means that Eqn.~\ref{errcontrol} cannot be uniform\footnote{We thank Pavel Bleher for pointing out a mistake in a former version of the article, which we corrected by introducing this family of nested contours.} in $n$ and $k$, even when $\beta = 2$.

A posteriori, from Eqn.~\ref{errcontrol}, we can deduce by choosing rather $l(k_0;n_0) = 8(n_0 + k_0)$:
\beq
\p \delta_{k_0} W_{n_0} \p \leq N^{-1}\,\p W_{n_0}^{\{k_0 + 1\}}\p_{\Gamma_{l(k_0 + 1;n_0)}} + \Delta_{n_0}^{\{k_0\}}\,\mathrm{max}(N^{-1}\,;\,\epsilon_{N}^{\{k_0\}}). \nn
\eeq
Subsequently, upon redefinition of the constant $\Delta_{n_0}^{\{k_0\}}$, we may choose $\epsilon_{N}^{\{k_0\}} = 1/N$. Finally, we can make the weaker statement that, for any $n$ and $k$:
\beq
\p \delta_{k} W_n \p _{\Gamma_E} \in o(1/N), \nn
\eeq
without uniformity in $n$ and $k$.

\section{Proof of the main results}
\label{final}

\subsection{Expansion of the correlators}
\label{fina}
We wish to study the $\beta$ ensembles on a given interval $[b_-,b_+]$, with the hypotheses~\ref{hypoge} on the potential $V$. When both edges are hard, Hyp.~\ref{hypoge} are equivalent to the five assumptions of Section~\ref{sec:rec1}, so the Proposition~\ref{prop1} is  already proved, as we have shown recursively that Eqn.~\ref{eq:expk0} holds for all $k_0$. Let us now assume that one of the edge is soft. The equilibrium measure $\mu_{\mathrm{eq}} := \mu_{\mathrm{eq}}^{V;[b_-,b_+]}$ with support $[\alpha_-,\alpha_+] \subset [b_-,b_+]$ also coincides with $\mu_{\mathrm{eq}}^{V;[a_-,a_+]}$, where $a_-$ can be any point in $[b_-,\alpha_-[$ if $b_-$ is a soft edge, and $a_- = b_-$ else (resp. $a_+$ can be any point in $]\alpha_+,b_+]$ if $b_+$ is a soft edge, and $a_+ = b_+$ else). When $b_{\tau}$ is a soft edge, "offcriticality" implies that $S(x)$ is positive in a neighborhood of $\alpha_{\tau}$ in $[b_-,b_+]\setminus]\alpha_-,\alpha_+[$. So, one can choose an interval $[a_-,a_+] \subseteq U$, and such that the five assumptions of Section~\ref{sec:rec1} are satisfied for $\dd\mu_{N,\beta}^{V;[a_-,a_+]}$. Theorem~\ref{recerror} then can be applied: there exists an asymptotic expansion
\beq
\label{fqy}W_n^{V;[a_-,a_+]}(x_1,\ldots,x_n) = \sum_{k \geq n - 2} N^{-k}\,W_n^{V;\{k\}}(x_1,\ldots,x_n),
\eeq
with respect to the norm $\p \cdot \p_{\Gamma_E}$ where $\Gamma_E \subseteq U$ can be any contour surrounding $[a_-,a_+]$ but not the zeroes of $S$. The "large deviation control" on $[b_-,b_+]$ allows to use Proposition~\ref{uuu2}: there exists $\eta > 0$ such that, for any contour $\Gamma_E' \subseteq \mathbb{C}$ surrounding $[b_-,b_+]$, there exists $T_{n,\Gamma} > 0$ such that:
\beq
\p W_n^{V;[b_-,b_+]} - W_n^{V;[a_-,a_+]} \p_{\Gamma_E'} \leq T_{n,\Gamma_E'}\,e^{-N\,\eta}. \nn
\eeq
This implies that the right hand side of Eqn.~\ref{fqy} is an asymptotic series for $W_n^{V;[b_-,b_+]}(x_1,\ldots,x_n)$, uniformly for $(x_1,\ldots,x_n)$ in any compact of $(\mathbb{C}\setminus[b_-,b_+])^n$.

We give below a more transparent condition, which imply the "large deviation control" assumption on $[b_-,b_+]$:
\begin{remark}
If $S(x) > 0$ whenever $x \in [b_-,b_+]$, then $\mathcal{J}^{V;[b_-,b_+]}$  achieves its minimum value only on $[\alpha_-,\alpha_+]$,
\end{remark}
Indeed, $\mathcal{J}^{V;[b_-,b_+]}(x)$ is differentiable when $x \in ]b_-,b_+[\setminus[\alpha_-,\alpha_+]$, and we have:
\beq
\big(\mathcal{J}^{V;[b_-,b_+]}(x)\big)' = \frac{(V^{\{0\}})'(x)}{2} - W_1^{\{-1\}}(x) = y(x) = S(x)\sigma(x). \nn
\eeq
The sign of the square root $\sigma(x)$ is determined for example by the positivity conditions \ref{ina} on $\mathcal{J}^{V;[b_-,b_+]}$. If we assume that $S$ do not vanish on $[b_-,b_+]$, this implies that  $\mathcal{J}^{V;[b_-,b_+]}$ is strictly decreasing in $[b_-,\alpha_-[$ and strictly increasing on $]\alpha_+,b_+]$, hence the remark.

\subsection{Expansion of the free energy}\label{freesec}

So far, we only have determined the expansion of the correlators which are by definition derivatives of the free energy. To find the free energy itself,
 one would like to interpolate  between our initial potential $V$, and a simpler situation, using  that the difference depends on the correlators.
For any fixed $\alpha_- < \alpha_+$, and fixed nature of the edges $\mathsf{X}_{\pm} \in \{\mathsf{hard},\mathsf{soft}\}$, we denote by $\mathcal{V}^{\alpha_+,\mathsf{X}_+}_{\alpha_-,\mathsf{X}_-}$ the set of potentials $V$:
 \begin{itemize}
 \item[$\bullet$] defined at least on some interval $[a_-,a_+] \supseteq [\alpha_-,\alpha_+]$, with $a_{\tau} \neq \alpha_{\tau}$ if $\mathsf{X}_{\tau} = \mathsf{soft}$, and $a_{\tau} = \alpha_{\tau}$ if $X_{\tau} = \mathsf{hard}$ ;
 \item[$\bullet$] which satisfies the five assumptions of Section~\ref{hypotheses} on $[a_-,a_+]$, in particular is offcritical on $[a_-,a_+]$ ;
 \item[$\bullet$] for which the equilibrium measure $\mu_{\mathrm{eq}}^{V;[a_-,a_+]}$ has $[\alpha_-,\alpha_+]$ as support,
 \item[$\bullet$] and such that $a_{\tau}$ is an edge of nature $\mathsf{X}_{\tau}$.
\end{itemize}

\begin{lemma}
$\mathcal{V}^{\alpha_+,\mathsf{X}_+}_{\alpha_-,\mathsf{X}_-}$ is a convex set.
\end{lemma}
\textbf{Proof.} Let $V_0, V_1 \in \mathcal{V}^{\alpha_+,\mathsf{X}_+}_{\alpha_-,\mathsf{X}_-}$, and set $V_s = (1 - s)V_0 + sV_1$ for $s \in [0,1]$. $V_0$ and $V_1$ are at least defined on a common interval $[a_-,a_+] \supseteq [\alpha_-,\alpha_+]$. Let us call $\nu_s = \dd\mu_{\mathrm{eq}}^{V_s;[a_-,a_+]}$ the equilibrium measure for the potential $V_s$ on $[a_-,a_+]$. We observe that $(1 - s)\dd \nu_0 + s\dd \nu_1$ is a probability measure which is solution of the characterization of $\dd L_s$ by Thm.~\ref{th:1}. Therefore, $\dd L_s = (1 - s)\dd L_0 + s \dd L_1$. Besides, we know that there exists a function $S_s$, regular in a neighborhood of $[\alpha_-,\alpha_+]$ in the complex plane, positive on $[a_-,a_+]$, such that:
\beq
\dd \nu_s(\xi) = \frac{\dd \xi}{\pi}\,S_s(\xi)\,\sqrt{\frac{\prod_{\tau\,/\,X_{\tau} = \mathsf{soft}} |\xi - \alpha_{\tau}|}{\prod_{\tau'\,/\,X_{\tau'} = \mathsf{hard}}|\xi - \alpha_{\tau'}|}}\,\mathbf{1}_{[\alpha_-,\alpha_+]}(\xi), \nn
\eeq
for $s = 0$ or $s = 1$.  Since the edges are of the same nature in $V_0$ et $V_1$, we must have $S_s = (1 - s)S_0 + sS_1$. Since $S_0$ and $S_1$ are positive on $[a_-,a_+]$, so is $S_s$. Hence $V_s \in \mathcal{V}^{\alpha_+,\mathsf{X}_+}_{\alpha_-,\mathsf{X}_-}$.

\begin{corollary}
\label{interpol}Let $V_0, V_1 \in \mathcal{V}^{\alpha_+,\mathsf{X}_+}_{\alpha_-,\mathsf{X}_-}$. When $a_-$ and $a_+$ satisfy the condition above, the quantity:
\beq
\ln Z_{N,\beta}^{V_1;[a_-,a_+]} - \ln Z_{N,\beta}^{V_0;[a_-,a_+]} = -\frac{N\beta}{2}\,\int_{0}^{1}\dd s\,\oint_{\mathcal{C}([a_-,a_+])} \frac{\dd\xi}{2i\pi}\,\big(V_1(\xi) - V_0(\xi)\big) W_1^{V_s;[a_-,a_+]}(\xi) \nn
\eeq
has a large $N$ asymptotic expansion of the form:
\beq
\ln Z_{N,\beta}^{V_1} - \ln Z_{N,\beta}^{V_0} = \sum_{k \geq {-2}} N^{-k}\,F^{V_0\rightarrow V_1;[a_-,a_+];\{k\}}_{\beta}, \nn
\eeq
where:
\beq
F^{V_0 \rightarrow V_1;[a_-,a_+];\{k\}}_{\beta} = -\frac{\beta}{2}\,\int_{0}^{1}\dd s \oint_{\mathcal{C}([a_-,a_+])} \frac{\dd\xi}{2i\pi}\sum_{m = 0}^{k + 2} \big(V_1^{\{m\}}(\xi) - V_0^{\{m\}}(\xi)\big)\,\big(W_1^{V_s;[a_-,a_+]}\big)^{\{k + 1 - m\}}(\xi). \nn
\eeq
\end{corollary}
\noindent\textbf{Proof.} Since $V_s$ satisfies the five assumptions of Section~\ref{hypotheses} for any $s \in [0,1]$, we can apply our main theorem to $W_1^{V_s}$. Moreover, since we do not reach a critical point when $s$ is in the compact $[0,1]$, we know that the error $O(N^{-K})$ made if we replace $W_1^{V_s}$ by $\sum_{k = -1}^{K - 1} N^{-k}\,W_1^{V_s;\{k\}}$ is uniformly bounded with respect to $s$ on some contour surrounding $[a_-,a_+]$ and in the analyticity domain of $V$. Therefore, we can exchange the integral and the sum in the asymptotic expansion.\hfill $\diamond$

\vspace{0.2cm}

For instance, when $V$ satisfies the five assumptions of Section~\ref{hypotheses} on some interval $[a_-,a_+]$, such that $a_{\pm}$ are soft edges, one can interpolate between $V$ and a Gaussian potential corresponding to an equilibrium measure with  support $[\alpha_-,\alpha_+]$:
\beq
V_{\mathrm{G},\alpha_-,\alpha_+}(x) = \frac{8}{(\alpha_+ - \alpha_-)^2}\Big(x - \frac{\alpha_- + \alpha_+}{2}\Big)^2. \nn
\eeq

\begin{proposition}
\label{maintho}Let $V$ be a potential satisfying the five assumptions of Section~\ref{hypotheses} on some interval $[a_-,a_+]$, such that $a_{\pm}$ are soft edges. For all $s \in [0,1]$, $(1 - s)V + sV_{\mathrm{G},\alpha_-,\alpha_+}$ belongs to $\mathcal{V}^{\alpha_+,\mathsf{soft}}_{\alpha_-,\mathsf{soft}}$ and we have the following asymptotic expansion when $N \rightarrow \infty$:
\beq
Z_{N,\beta}^{V} = Z_{N,\beta}^{V_{\mathrm{G},\alpha_-,\alpha_+}}\,\exp\Big(\sum_{k \geq -2} N^{-k}\,F^{V\rightarrow V_{\mathrm{G},\alpha_-,\alpha_+};[a_-,a_+];\{k\}}_{\beta}\Big), \nn
\eeq
where the prefactor is a partition function of the Gaussian $\beta$ ensemble (see Eqn.~\ref{eq:Selb}):
\beq
Z_{N,\beta}^{V_{\mathrm{G},\alpha_-,\alpha_+}} = Z_{N,\mathrm{G}\beta\mathrm{E}}\,\Big(\frac{\alpha_+ - \alpha_-}{4}\Big)^{N + \frac{\beta}{2}\,N(N - 1)}. \nn
\eeq
\end{proposition}
According to the discussion of \S~\ref{fina}, we can weaken the hypothesis of the proposition above to find Theorem~\ref{prop0}.

\subsection{Central limit theorem}\label{secclt}
Eventually, our results imply the central limit theorem proved by Johansson \cite{Johansson}, but
here integration is taken on a compact set $[a_-,a_+]$ instead of the
real line (in fact as our derivation is quite similar to Johansson's, this is not surprising).
For simplicity, we take here the hypotheses of Section~\ref{sec:rec1}, although we could refine to hypotheses~\ref{hypoge} following \S~\ref{fina}.

Let $h\,:\, [a_-,a_+] \rightarrow \mathbb{R}$ be a function which can be extended as a holomorphic function defined on some neighborhood of $[a_-,a_+]$, let us take $V \equiv V^{\{0\}}$ independent of $N$, $V^{\{1\}} = \frac{2}{\beta}\:h$ and define $V_h = V^{\{0\}} + N^{-1}\,V^{\{1\}} = V -\frac{2}{N\beta}\:h$. Then:
\beq
\mu^{V;[a_-,a_+]}_{N,\beta}\left[\exp\Big(\sum_{i=1}^N h(\lambda_i)\Big)\right] =
\frac{Z_{N,\beta}^{V_h;[a_-,a_+]}}{Z_{N,\beta}^{V;[a_-,a_+]}}, \nn
\eeq
and we can use Corollary~\ref{interpol} to derive its large $N$ asymptotics.
Indeed, we have
\beq
\ln \mu^{V;[a_-,a_+]}_{N,\beta}\left[\exp\Big(\sum_{i=1}^N h(\lambda_i)\Big)\right] = \int_0^1 \dd s\,\oint_{{\mathcal C }([a_-,a_+])}\frac{\dd\xi}{2i\pi}\,W_{1}^{V_{sh}}(\xi)\, h(\xi). \nn
\eeq
By Theorem \ref{maintho}, or simply at the point of Lemma~\ref{lema}, we have:
\bea
 W_1^{V_{sh};\{-1\}} (\xi) & = & W_{1}^{V;\{-1\}}(\xi) = \int \frac{\dd \mu_{\mathrm{eq}}(\eta)}{\xi - \eta}, \nn \\
 W_1^{V_{sh};\{0\}} (\xi) & = & \mathcal{K}^{-1}\left\{-\Big(1 - \frac{2}{\beta}\Big)\left[\frac{\dd}{\dd x}\big(W_1^{V;\{-1\}}(x)\big) + \sum_{\tau \in \mathsf{Hard}} \frac{1}{a_\tau - a_{-\tau}}\,\frac{1}{x - a_{\tau}}\right] \right. \nn \\
 && \left. \phantom{\mathcal{K}^{-1}\,} - \frac{2s}{\beta}\mathcal{N}_{h'}(W_1^{V;\{-1\}})(x)\right\}, \nn \\
 W_{1}^{V_{sh}} & = & N\,W_1^{V_{sh};\{-1\}} + W_{1}^{V_{sh};\{0\}} + o(1), \nonumber \eea
which shows the:
\begin{proposition}
\label{propclt}Central limit theorem.
{\small\beq\ln \mu^{V;[a_-,a_+]}_{N,\beta}\left[\exp\Big(\sum_{i=1}^N h(\lambda_i)\Big)\right] =
N\int \dd\mu_{\mathrm{eq}}(\eta)\,h(\eta) + m[h]+ \frac{1}{2}\,C[h]+o(1), \nn \eeq}
$\!\!\!$with $m[h]$ the linear in the function $h$, given by:
\beq
m[h]=-\Big(1 - \frac{2}{\beta}\Big)\oint_{{\mathcal C }(
[a_-,a_+])}\frac{\dd\xi}{2i\pi}\,\mathcal{K}^{-1}\left\{\frac{\dd}{\dd x}\big(W_1^{V;\{-1\}}(x)\big) + \sum_{\tau \in \mathsf{Hard}} \frac{1}{a_\tau - a_{-\tau}}\,\frac{1}{x - a_{\tau}}\right\}\,h(\xi), \nn
\eeq
and $C[h]$ the quadratic function of $h$ given by:
\beq
C[h]=-\frac{2}{\beta} \oint_{{\mathcal C }([a_-,a_+])}\frac{\dd\xi}{2i\pi}\,\mathcal{K}^{-1}\Big[
 \mathcal{N}_{h'}(W_{1}^{V;\{-1\}})\Big](\xi)\,h(\xi). \nn
\eeq
Therefore $\sum_{i=1}^N h(\lambda_i)-N\int \dd\mu_{\mathrm{eq}}(\eta)\,h(\eta)$
converges towards a Gaussian variable with mean $m[h]$ and covariance $C[h]$.
\end{proposition}

\appendix
\section{Proof of Proposition \ref{uuu3}}

We use the notation introduced in Proposition~\ref{uuu3}, in particular the eigenvalues are integrated over a segment $[b_-,b_+]$ which may not be compact.

\subsection{$\widetilde{\mathcal{J}}^{V;[b_-,b_+]}$ is a good rate function}

$\widetilde{\mathcal{J}}^{V;[b_-,b_+]}$ is lower semicontinuous as a supremum of the continuous functions
$$\widetilde{\mathcal{J}}^{V;[b_-,b_+]}_\varepsilon(x):=\frac{V(x)}{2} - \int_{b_-}^{b_+} \dd \mu_{\mathrm{eq}}^{V;[b_-,b_+]}(\xi)\,\ln\big[\mathrm{max}(|x - \xi|,\varepsilon)\big] - \inf_{\xi \in [b_-,b_+]} \mathcal{J}^{V;[b_-,b_+]}(\xi).$$
Moreover, by the assumption of Eqn.~\ref{assss}, it goes to infinity at infinity. Hence, $\widetilde{\mathcal{J}}^{V;[b_-,b_+]}$ has compact level sets. Since it is non-negative, it is a good rate function.

\subsection{The law of the extreme eigenvalues is exponentially tight}

Exponential tightness of the extreme eigenvalues means:
\beq\label{expt}\limsup_{M\ra\infty}\limsup_{N\ra\infty} \frac{1}{N}\ln \mu^{V;[b_-,b_+]}_{N,\beta}\left( \lambda_{\mathrm{max}}\ge M\mbox{ or }
\lambda_{\mathrm{min}}\le - M\right)=-\infty.\eeq
By \cite[Lemma 2.6.7]{AGZ}, it is enough to show that:
\beq\label{res1}\limsup_{N\ra\infty} \frac{1}{N}\ln
\frac{Z^{V;[b_-,b_+]}_{N-1,\beta}}{Z^{V;[b_-,b_+]}_{N,\beta}}<\infty.
\eeq
For this purpose, observe that by Jensen's inequality
{\small\bea
\frac{Z^{V;[b_-,b_+]}_{N,\beta}}{Z^{V;[b_-,b_+]}_{N-1,\beta}}&=&\mu^{V;[b_-,b_+]}_{N-1,\beta}\left[\int_{b_-}^{b_+} \dd \lambda_N \exp\Big(\beta\sum_{i=1}^{N - 1} \ln|\lambda_N-\lambda_i|
-\frac{\beta N}{2}V(\lambda_N)-\frac{\beta}{2}\sum_{i=1}^{N-1}V(\lambda_i)\Big)\right] \nn \\
& \geq & \kappa\,\exp\left\{\frac{\beta}{2}\big(\mu_{N - 1,\beta}^{V;[b_-,b_+]}\otimes\chi\big)\left[2\sum_{i = 1}^{N - 1} \ln|\lambda_N - \lambda_i| - (N - 1)V(\lambda_N) - \sum_{i = 1}^{N - 1} V(\lambda_i)\right]\right\}, \nn
\eea}
$\!\!\!$where we denoted $\chi$ the law on $\lambda_N$ given by:
\beq
\dd\chi(x) = \frac{\mathbf{1} _{[b_-,b_+]}(x)\dd x}{\kappa}\,e^{-\frac{\beta}{2}V(x)},\qquad \kappa = \int_{b_-}^{b_+} \dd\xi\,e^{-\frac{\beta}{2} V(\xi)}. \nn
\eeq
The function $\xi \mapsto \int_{\mathbb{R}} \dd\chi(\lambda_N)\,\ln|\lambda_N - \xi|$ is bounded on compact sets and going to infinity like $\ln |\xi|$, so is bounded from below, by a constant $\frac{\kappa_1}{2}$. We can rewrite:
\beq
\frac{Z^{V;[b_-,b_+]}_{N,\beta}}{Z^{V;[b_-,b_+]}_{N -1,\beta}} \geq \kappa\,\exp\left\{\beta(N - 1)\Big[\kappa_1 - \chi[V] - \mu_{N - 1,\beta}^{V;[b_-,b_+]}[L_{N - 1}(V)]\Big]\right\}. \nn
\eeq
By exponential tightness \cite[Eqn. 2.6.21]{AGZ}, we know that there exists  a constant $\kappa_2> 0$ so that
\bea
-\mu_{N - 1,\beta}^{V;[b_-,b_+]}\big[L_{N-1}(V)\big]  & \geq & - \mu_{N - 1,\beta}^{V; [b_-,b_+]}\big[L_{N - 1}(|V|)\big] \geq  -\kappa_2. \nn
\eea
So, if we set $\kappa_3 = \chi[V]$ and choose $\kappa_2$ large enough, we have:
\beq
\frac{Z^{V;[b_-,b_+]}_{N,\beta}}{Z^{V; [b_-,b_+]}_{N - 1,\beta}} \geq \kappa\,e^{-\beta(N - 1)\delta}, \nn
\eeq
with a positive constant $\delta = -\kappa_1 + \kappa_2 + \kappa_3$. This justifies Eqn.~\ref{res1} and completes the proof of
Eqn.~\ref{expt}.

\subsection{Upper bound for large deviation of the extreme \mbox{eigenvalues}}

We give the argument for the minimal eigenvalue, the case of the maximal eigenvalue being similar. By exponential tightness (Eqn.~\ref{expt}),
it is enough to prove a weak large deviation upper bound, that is control
the probability of small balls.  First, observe that for any $x-\alpha_-\ge 2\epsilon>0$,

$$\mu_{N,\beta}^{V;[b_-,b_+]}[\lambda_{\mathrm{min}} \ge x]\le \mu_{N,\beta}^{V;[b_-,b_+]}[L_N(\mathbf{1}_{[\alpha_-,\alpha_-+\e]})=0]$$
is of order $e^{-N^2 \kappa_\epsilon}$ for some $\kappa_\epsilon>0$
by the large deviation principle for the law of $L_N$ under $ \mu_{N,\beta}^{V;[b_-,b_+]}$, see e.g. \cite{BAG} or  \cite[Theorem 2.6.1]{AGZ}. Moreover, the probability
that $\lambda_{\rm min}$ is smaller than $a_-$ vanishes and therefore we have
$$\limsup_{\epsilon\downarrow 0}
\limsup_{N\ra\infty}\frac{1}{N}\ln\mu^{V;[b_-,b_+]}_{N,\beta}\left(\lambda_{\rm min}\in ]-\infty,b_--\e]\cup [\alpha_-+\epsilon,+\infty[\right)=-\infty.$$
Hence, we may and shall concentrate on probability of deviating on $[b_-,\alpha_-]$,
and actually we may restrict ourselves to the case where
$b_-$ and $b_+$  are finite by Eqn.~\ref{expt}. We let $F$ be a closed subset of $[b_-,\alpha_-]$. We then have:
 \beq
\mu_{N,\beta}^{V;[b_-,b_+]}[\lambda_{\mathrm{min}} \in F]=
  Y_{N} \int_F
\dd\xi\,e^{-\frac{\beta}{2} V(\xi)}\,\Xi_{N}(\xi), \nn
\eeq
where we introduced:
\bea
 Y_{N} &=& \frac{Z^{\frac{N}{N-1}V;[b_-,b_+]}_{N-1,\beta}}{Z^{V;[b_-,b_+]}_{N,\beta}},\nn\\
\Xi_{N}(\xi)&=& \mu_{N-1,\beta}^{\frac{N}{N-1}V;[b_-,b_+]}\left(e^{\beta\sum_{i=1}^{N-1} \ln|\xi-\lambda_i|-\frac{\beta}{2}(N-1)V(\xi)}\,\prod_{i=1}^{N-1} \mathbf{1}_{[b_-,\lambda_i]}(\xi)\right).\nn
\eea

\subsubsection*{Upper bound for $\Xi_{N}(\xi)$}
Notice that the logarithm is uniformly bounded from above on compacts
so that the exponent is at most of order $N$. Therefore, we may and shall assume that under
$\mu _{N-1,\beta}^{\frac{N}{N-1}V;[b_-,b_+]}$, $L_{N-1}$ is at a distance smaller than $\kappa>0$ from the equilibrium measure $\mu_{\mathrm{eq}} := \mu_{\mathrm{eq}}^{V;[b_-,b_+]}$,
since the opposite event has probability smaller than $e^{-\Gamma_\kappa (N-1)^2}$ for some $\Gamma_\kappa >0$, see e.g. \cite[Theorem 2.6.1]{AGZ}. Here, the distance can be taken to be any distance compatible with the weak topology, e.g. the Wasserstein distance. Thus, we have for large $N$:
 \beq
\Xi_{N}(\xi)
 \le e^{-\Gamma_\kappa N^2/2}
 + e^{\beta(N-1)\sup_{d(\mu,\mu_{\mathrm{eq}})<\kappa} \big(-\frac{V(\xi)}{2} + \int\ln |\xi-\eta|\dd\mu(\eta)\big)}, \nn
\eeq
where we take the supremum over probability measures $\mu$ on $[b_-,b_+]$ with Wasserstein distance to $\mu_{eq}$ strictly smaller than $\kappa$. We observe also that for all probability measures
$\mu$ on $[b_-,b_+]$, and for any $\zeta>0$:
\beq
\int_{b_-}^{b_+}\ln |\xi-\eta|\dd\mu(\eta)\le \phi_\zeta(\mu,\xi)
=\int_{b_-}^{b_+} \ln\big[\mathrm{max}(|\xi-\eta|,\zeta)\big]\dd\mu(\eta) \,\nn
\eeq
 where $\phi_\zeta(\mu,\xi)$
is  continuous in $\mu$ and $\xi$,
and $\phi_\zeta(\mu_{\mathrm{eq}},\xi)$ converges towards $\phi_0(\mu_{\mathrm{eq}},\xi)$ as $\zeta$ goes to zero. We deduce that:
\beq
\limsup_{\kappa\downarrow 0}\sup_{\xi\in F}\sup_{d(\mu,\mu_{\mathrm{eq}})<\kappa} \beta\left(\int\ln |\xi-\eta|d\mu(\eta) -\frac{V(\xi)}{2}\right)\le
-\beta\,\inf_{\xi \in F} \mathcal{J}^{V;[b_-,b_+]}(\xi). \nn
\eeq
Therefore, for any $\eta'>0$, and $N$ large enough, we conclude that:
\beq\label{cont2}
\sup_{\xi\in F}  \Xi_N(\xi )\le e^{N\big(\eta' - \beta\inf_{\xi \in F} \mathcal{J}^{V;[b_-,b_+]}(\xi)\big)}.
\eeq

\subsubsection*{Lower bound for $Y_N$}

We observe that, for any $\varepsilon>0$ small enough, and any $x \in [b_- + \varepsilon,b_+ - \varepsilon]$, there exists $\delta_{\varepsilon}$ going to zero with $\varepsilon$
so that
\bea
\frac{1}{Y_N} &=&\frac{Z^{V;[b_-,b_+]}_{N,\beta}}{Z^{\frac{N}{N-1}V;[b_-,b_+]}_{N-1,\beta}} \nn \\
&=& \mu^{\frac{N}{N-1}V;[b_-,b_+]}_{N-1,\beta}\left(\int_{b_-}^{b_+} \dd\xi\,e^{-\frac{\beta N}{2} V(\xi) }\prod_{i=1}^{N-1}|\xi-\lambda_i|^\beta\right)\nn\\
&\ge&  \mu^{\frac{N}{N-1}V;[b_-,b_+]}_{N-1,\beta}\left(\int_{x-\varepsilon}^{x+\varepsilon}\dd\xi\,e^{-\frac{\beta N}{2} V(\xi) }\prod_{i=1}^{N-1}|\xi-\lambda_i|^\beta\right)\nn\\
&\ge& 2\varepsilon\,e^{-\frac{\beta N}{2} V(x) -N\delta_\varepsilon}\mu^{\frac{N}{N-1}V;[b_-,b_+]}_{N-1,\beta}\left(e^{\sum_{i=1}^{N-1}\frac{\beta}{2\varepsilon}\int_{x-\varepsilon}^{x+\varepsilon}
\ln|\xi-\lambda_i|\,\dd\xi}\right),\nn
\eea
where we have finally used Jensen's inequality. But $\lambda\rightarrow \frac{1}{2\varepsilon}\int_{x-\varepsilon}^{x+\varepsilon}
\ln|\xi-\lambda| d\xi $ is bounded continuous on $[a_-,a_+]$ and therefore by the large deviation principle for the law of the empirical measure
$L_{N-1}$ under $\mu^{\frac{N}{N-1}V;[b_-,b_+]}_{N-1,\beta}$  (with rate function which vanishes only at $\mu_{\mathrm{eq}}$) we deduce that for $N$ large enough:
\beq
\frac{1}{Y_N} \ge  2\varepsilon\,e^{-\frac{\beta N}{2} V(x) -2N\delta_\varepsilon}\,e^{(N-1)\int \frac{\beta}{2\varepsilon}\big(\int_{x-\varepsilon}^{x+\varepsilon}
\ln|\xi-\lambda|\dd\xi\big)\dd\mu_{\mathrm{eq}}(\lambda)}.\nn
\eeq
Hence, by taking $\varepsilon$ sufficiently small independently of $N$, and optimizing over the choice of $x \in ]b_-,b_+[$, we conclude that for any $\eta''>0$, and $N$ large enough,
\beq\label{cont3}
\frac{1}{Y_N} \ge e^{-N\big(\eta'' + \beta \inf_{\xi \in [b_-,b_+]} \mathcal{J}^{V;[b_-,b_+]}(\xi)\big)}.
\eeq
Putting Eqn.~\ref{cont2} and \ref{cont3} together, we deduce that for all $\delta>0$ and $N$ large enough:
\beq
\mu^{V;[b_-,b_+]}_{N,\beta}\left(\lambda_{\rm min}\in F\right)\le
e^{N\beta\big(-\inf_{x \in F} \mathcal{J}^{V;[b_-,b_+]}(x) + \inf_{\xi \in [b_-,b_+]} \mathcal{J}^{V;[b_-,b_+]}(\xi) + \delta\big)}, \nn
\eeq
which provides the announced upper bound.

\subsubsection*{Conclusion}

As a consequence, since we assumed that the rate function only vanishes at $\alpha_-,\alpha_+$ we
deduce that for any $\epsilon>0$, there exists $\delta_\epsilon>0$ so that:
\beq\label{expb}
\mu^{V,[b_-,b_+]}_{N,\beta}\left(\lambda_{\rm min}\le \alpha_--\epsilon\right)\le e^{-\delta_\epsilon N},
\eeq
as well as a similar result for the largest eigenvalue.

\subsection{Lower bound for large deviation of extreme eigenvalues}

To establish a lower bound, we start again from Eqn.~\ref{theeq} with an open ball $B= ]x - \epsilon,x + \epsilon[ \subset [b_-,\alpha_-]$:
\beq
\mu_{N,\beta}^{V,[b_-,b_+]}\left( \lambda_{\rm min}\in B\right)
=Y_N\int _B \dd\xi e^{-\frac{\beta}{2}V(\xi)}\Xi_N(\xi), \nn
\eeq
but replace the role of $Y_N$ and $\Xi_N$ in the bounds.
Namely,
we first have by Jensen's inequality:
\beq
\int _B \dd\xi\,e^{-\frac{\beta}{2}V(\xi)}\Xi_N(\xi)\ge \kappa_N
e^{ \int  \dd\widetilde{\chi}(\xi,\lambda)\big(
\beta\sum_{i=1}^{N-1} \ln|\xi-\lambda_i|-\frac{\beta}{2}(N-1)V(\xi)\big)}, \nn\eeq
with
\bea
\dd\widetilde{\chi}(\xi,\lambda) & = & \frac{\mathbf{1}_B(\xi)\,\mathbf{1}_{\lambda_{\mathrm{min}}\ge \xi}}{\kappa_N}
 \dd\xi\,e^{-\frac{\beta}{2}V(\xi)}\,\dd\mu_{N-1,\beta}^{\frac{N}{N-1}V;[b_-,b_+]}(\lambda), \nn \\
 \kappa_N & = &
\int_B \dd\xi\,e^{-\frac{\beta}{2}V(\xi)}
 \mu_{N-1,\beta}^{\frac{N}{N-1}V;[b_-,b_+]}[\mathbf{1}_{\lambda_{\mathrm{min}} \ge \xi}]. \nn
\eea
Thanks to Eqn.~\ref{expb} (note that it applies similarly to $NV/(N-1)$ as the assumptions
does not depend on the fine asymptotics of $V$), we know that $\kappa_N$ converges
towards a non vanishing constant. Moreover, the logarithm, once integrated
against $\dd\xi$, produces a smooth bounded function and therefore we can use the convergence of $L_{N-1}$
towards $\mu_{\rm eq}$ under $\mu_{N-1,\beta}^{\frac{N}{N-1}V;[b_-,b_+]}$ to conclude that:
\beq
\liminf_{N\ra\infty}\frac{1}{N}\ln
\int _B\dd\xi\,e^{-\frac{\beta}{2}V(\xi)}\Xi_N(\xi)\ge-\frac{\beta}{2}
\,\frac{ \int_B \dd\xi\, e^{-\frac{\beta}{2}V(\xi)}
\Big(V(\xi)-2 \int \dd\mu_{\rm eq}(\eta)\ln |\xi-\eta|\Big)}{\int_B \dd\xi\,e^{-\frac{\beta}{2}V(\xi)} }.\nn\eeq
Letting now $\e$ going to zero in
 $B= ]x - \e,x + \e[$ proves that:
\beq\label{poi}
\liminf_{\epsilon\ra 0}\liminf_{N\ra\infty}\frac{1}{N}\ln
\int _{B(x,\epsilon)} \dd\xi\,e^{-\frac{\beta}{2}V(\xi)}\,\Xi_N(\xi)\ge  \beta\Big( \int\dd\mu_{\rm eq}(\eta)\,\ln |\xi-\eta| -\frac{V(\eta)}{2} \Big). \eeq
To bound $Y_N$ from below, it is enough to bound $1/Y_N$ from above, which can be done in the same way we bounded $\Xi_N$ from above in the argument for the upper bound.
We finally conclude:
\beq
 \liminf_{\epsilon\ra 0}\liminf_{N\ra\infty}\frac{1}{N}\ln
\mu^{V,[b_-,b_+]}_{N,\beta}\left(\lambda_{\rm min}\in ]x-\epsilon,x + \epsilon[\right)
\ge -\beta\widetilde{\mathcal{J}}^{V;[b_-,b_+]}(x), \nn\eeq
which completes the proof of the large deviation principle.

\section*{Acknowledgments}
 We would like to thank the MSRI and the organizers of the semester "Random Matrix Theory and its Applications" where this work was initiated, as well as Bertrand Eynard and Pavel Bleher for fruitful discussions. This work was supported by the ANR project ANR-08-BLAN-0311-01. The work of G.B. is supported by the SWISS NSF (no $200021_143434$) and the ERC AG CONFRA.

\newpage

\bibliographystyle{amsalpha}
\bibliography{Bib240411}

\end{document}